\documentclass[12pt]{article}

\catcode`\@=11

\thicklines
\newskip\Einheit \Einheit=.6cm
\newcount\xcoord \newcount\ycoord
\newdimen\xdim \newdimen\ydim \newdimen\PfadD@cke \newdimen\Pfadd@cke
\def\PfadDicke#1{\PfadD@cke#1 \divide\PfadD@cke by2 
\Pfadd@cke\PfadD@cke \multiply\PfadD@cke by2}
\long\def\LOOP#1\REPEAT{\def\BODY{#1}\ITERATE}
\def\ITERATE{\BODY \let\next\ITERATE \else\let\next\relax\fi \next}
\let\REPEAT=\fi
\def\Punkt{\hbox{\raise-2pt\hbox to0pt{\hss\scriptsize$\bullet$\hss}}}

\def\DuennPunkt(#1,#2){\unskip
  \raise#2 \Einheit\hbox to0pt{\hskip#1 \Einheit
          \raise-1.5pt\hbox to0pt{\hss\tiny$\bullet$\hss}\hss}}
		  
\def\NormalPunkt(#1,#2){\unskip
  \raise#2 \Einheit\hbox to0pt{\hskip#1 \Einheit
          \raise-3pt\hbox to0pt{\hss\large$\bullet$\hss}\hss}}
\def\DickPunkt(#1,#2){\unskip
  \raise#2 \Einheit\hbox to0pt{\hskip#1 \Einheit
          \raise-4pt\hbox to0pt{\hss\Large$\bullet$\hss}\hss}}
\def\Kreis(#1,#2){\unskip
  \raise#2 \Einheit\hbox to0pt{\hskip#1 \Einheit
          \raise-4pt\hbox to0pt{\hss\Large$\circ$\hss}\hss}}
\def\Diagonale(#1,#2)#3{\unskip\leavevmode
  \xcoord#1\relax \ycoord#2\relax
      \raise\ycoord \Einheit\hbox to0pt{\hskip\xcoord \Einheit
         \unitlength\Einheit
         \line(1,1){#3}\hss}}
\def\AntiDiagonale(#1,#2)#3{\unskip\leavevmode
  \xcoord#1\relax \ycoord#2\relax \advance\xcoord by -0.05\relax
      \raise\ycoord \Einheit\hbox to0pt{\hskip\xcoord \Einheit
         \unitlength\Einheit
         \line(1,-1){#3}\hss}}
\def\Pfad(#1,#2),#3\endPfad{\unskip\leavevmode
  \xcoord#1 \ycoord#2 \thicklines\ZeichnePfad#3\endPfad\thinlines}
\def\ZeichnePfad#1{\ifx#1\endPfad\let\next\relax
  \else\let\next\ZeichnePfad
    \ifnum#1=1
      \raise\ycoord \Einheit\hbox to0pt{\hskip\xcoord \Einheit
         \vrule height\Pfadd@cke width1 \Einheit depth\Pfadd@cke\hss}%
      \advance\xcoord by 1
     \else\ifnum#1=2
      \raise\ycoord \Einheit\hbox to0pt{\hskip\xcoord \Einheit
         \unitlength\Einheit
         \line(0,1){1}\hss}
      \advance\xcoord by 0
      \advance\ycoord by 1
 \else\ifnum#1=3
      \raise\ycoord \Einheit\hbox to0pt{\hskip\xcoord \Einheit
         \unitlength\Einheit
         \line(1,1){1}\hss}
      \advance\xcoord by 1
      \advance\ycoord by 1
    \else\ifnum#1=4
      \raise\ycoord \Einheit\hbox to0pt{\hskip\xcoord \Einheit
         \unitlength\Einheit
         \line(1,-1){1}\hss}
      \advance\xcoord by 1
      \advance\ycoord by -1
   \else\ifnum#1=5
      \raise\ycoord \Einheit\hbox to0pt{\hskip\xcoord \Einheit
         \unitlength\Einheit
         \line(2,1){2}\hss}
      \advance\xcoord by 2
      \advance\ycoord by 1
	  \else\ifnum#1=6
      \raise\ycoord \Einheit\hbox to0pt{\hskip\xcoord \Einheit
         \unitlength\Einheit
         \line(2,-1){2}\hss}
      \advance\xcoord by 2
      \advance\ycoord by -1
	  \else\ifnum#1=7
      \raise\ycoord \Einheit\hbox to0pt{\hskip\xcoord \Einheit
         \unitlength\Einheit
         \line(3,1){3}\hss}
      \advance\xcoord by 3
      \advance\ycoord by 1
	  \else\ifnum#1=8
      \raise\ycoord \Einheit\hbox to0pt{\hskip\xcoord \Einheit
         \unitlength\Einheit
         \line(3,-1){3}\hss}
      \advance\xcoord by 3
      \advance\ycoord by -1
    \fi\fi\fi\fi\fi\fi\fi\fi
  \fi\next}
\def\hSSchritt{\leavevmode\raise-.4pt\hbox 
to0pt{\hss.\hss}\hskip.2\Einheit
  \raise-.4pt\hbox to0pt{\hss.\hss}\hskip.2\Einheit
  \raise-.4pt\hbox to0pt{\hss.\hss}\hskip.2\Einheit
  \raise-.4pt\hbox to0pt{\hss.\hss}\hskip.2\Einheit
  \raise-.4pt\hbox to0pt{\hss.\hss}\hskip.2\Einheit}
\def\vSSchritt{\vbox{\baselineskip.2\Einheit\lineskiplimit0pt
\hbox{.}\hbox{.}\hbox{.}\hbox{.}\hbox{.}}}
\def\DSSchritt{\leavevmode\raise-.4pt\hbox to0pt{%
  \hbox to0pt{\hss.\hss}\hskip.2\Einheit
  \raise.2\Einheit\hbox to0pt{\hss.\hss}\hskip.2\Einheit
  \raise.4\Einheit\hbox to0pt{\hss.\hss}\hskip.2\Einheit
  \raise.6\Einheit\hbox to0pt{\hss.\hss}\hskip.2\Einheit
  \raise.8\Einheit\hbox to0pt{\hss.\hss}\hss}}
\def\dSSchritt{\leavevmode\raise-.4pt\hbox to0pt{%
  \hbox to0pt{\hss.\hss}\hskip.2\Einheit
  \raise-.2\Einheit\hbox to0pt{\hss.\hss}\hskip.2\Einheit
  \raise-.4\Einheit\hbox to0pt{\hss.\hss}\hskip.2\Einheit
  \raise-.6\Einheit\hbox to0pt{\hss.\hss}\hskip.2\Einheit
  \raise-.8\Einheit\hbox to0pt{\hss.\hss}\hss}}
\def\SPfad(#1,#2),#3\endSPfad{\unskip\leavevmode
  \xcoord#1 \ycoord#2 \ZeichneSPfad#3\endSPfad}
\def\ZeichneSPfad#1{\ifx#1\endSPfad\let\next\relax
  \else\let\next\ZeichneSPfad
    \ifnum#1=1
      \raise\ycoord \Einheit\hbox to0pt{\hskip\xcoord \Einheit
         \hSSchritt\hss}%
      \advance\xcoord by 1
    \else\ifnum#1=2
      \raise\ycoord \Einheit\hbox to0pt{\hskip\xcoord \Einheit
        \hbox{\hskip-2pt \vSSchritt}\hss}%
      \advance\ycoord by 1
    \else\ifnum#1=3
      \raise\ycoord \Einheit\hbox to0pt{\hskip\xcoord \Einheit
         \DSSchritt\hss}
      \advance\xcoord by 1
      \advance\ycoord by 1
    \else\ifnum#1=4
      \raise\ycoord \Einheit\hbox to0pt{\hskip\xcoord \Einheit
         \dSSchritt\hss}
      \advance\xcoord by 1
      \advance\ycoord by -1
    \fi\fi\fi\fi
  \fi\next}
\def\Koordinatenachsen(#1,#2){\unskip
 \hbox to0pt{\hskip-.5pt\vrule height#2 \Einheit width.5pt depth1 
\Einheit}%
 \hbox to0pt{\hskip-1 \Einheit \xcoord#1 \advance\xcoord by1
    \vrule height0.25pt width\xcoord \Einheit depth0.25pt\hss}}
\def\Koordinatenachsen(#1,#2)(#3,#4){\unskip
 \hbox to0pt{\hskip-.5pt \ycoord-#4 \advance\ycoord by1
    \vrule height#2 \Einheit width.5pt depth\ycoord \Einheit}%
 \hbox to0pt{\hskip-1 \Einheit \hskip#3\Einheit 
    \xcoord#1 \advance\xcoord by1 \advance\xcoord by-#3 
    \vrule height0.25pt width\xcoord \Einheit depth0.25pt\hss}}
\def\Gitter(#1,#2){\unskip \xcoord0 \ycoord0 \leavevmode
  \LOOP\ifnum\ycoord<#2
    \loop\ifnum\xcoord<#1
      \raise\ycoord \Einheit\hbox to0pt{\hskip\xcoord 
\Einheit\Punkt\hss}%
      \advance\xcoord by1
    \repeat
    \xcoord0
    \advance\ycoord by1
  \REPEAT}
\def\Gitter(#1,#2)(#3,#4){\unskip \xcoord#3 \ycoord#4 \leavevmode
  \LOOP\ifnum\ycoord<#2
    \loop\ifnum\xcoord<#1
      \raise\ycoord \Einheit\hbox to0pt{\hskip\xcoord 
\Einheit\Punkt\hss}%
      \advance\xcoord by1
    \repeat
    \xcoord#3
    \advance\ycoord by1
  \REPEAT}
\def\Label#1#2(#3,#4){\unskip \xdim#3 \Einheit \ydim#4 \Einheit
  \def\lo{\advance\xdim by-.5 \Einheit \advance\ydim by.5 \Einheit}%
  \def\llo{\advance\xdim by-.25cm \advance\ydim by.5 \Einheit}%
  \def\loo{\advance\xdim by-.5 \Einheit \advance\ydim by.25cm}%
  \def\o{\advance\ydim by.25cm}%
  \def\ro{\advance\xdim by.5 \Einheit \advance\ydim by.5 \Einheit}%
  \def\rro{\advance\xdim by.25cm \advance\ydim by.5 \Einheit}%
  \def\roo{\advance\xdim by.5 \Einheit \advance\ydim by.25cm}%
  \def\l{\advance\xdim by-.30cm}%
  \def\r{\advance\xdim by.30cm}%
  \def\lu{\advance\xdim by-.5 \Einheit \advance\ydim by-.6 \Einheit}%
  \def\llu{\advance\xdim by-.25cm \advance\ydim by-.6 \Einheit}%
  \def\luu{\advance\xdim by-.5 \Einheit \advance\ydim by-.30cm}%
  \def\u{\advance\ydim by-.30cm}%
  \def\ru{\advance\xdim by.5 \Einheit \advance\ydim by-.6 \Einheit}%
  \def\rru{\advance\xdim by.25cm \advance\ydim by-.6 \Einheit}%
  \def\ruu{\advance\xdim by.5 \Einheit \advance\ydim by-.30cm}%
  #1\raise\ydim\hbox to0pt{\hskip\xdim
     \vbox to0pt{\vss\hbox to0pt{\hss$#2$\hss}\vss}\hss}%
}
\catcode`\@=12

\usepackage{amsmath,amsthm}
\usepackage{amssymb}
\usepackage{color}
\usepackage{xspace}
\usepackage[colorlinks=true,
linkcolor=green,
filecolor=brown,
citecolor=green]{hyperref}
\def\red{\textcolor{red} }
\def\blue{\textcolor{blue} }
\def\green{\textcolor{green} }
\def\magenta{\textcolor{magenta} }

\def\ep{\epsilon}
\def\si{\sigma}
\def\DD{\ensuremath{\mathcal D}\xspace}

\def\gl{ground level}
\def\A{\ensuremath{\mathsf{A}}\xspace}
\def\D{\ensuremath{\mathsf{D}}\xspace}
\def\a{\ensuremath{\mathsf{a}}\xspace}
\def\d{\ensuremath{\mathsf{d}}\xspace}

\def\mbf#1{\mathchoice{\hbox{\boldmath $\displaystyle #1$}}
        {\hbox{\boldmath $\textstyle #1$}}
        {\hbox{\boldmath $\scriptstyle #1$}}
        {\hbox{\boldmath $\scriptscriptstyle #1$}}}

\begin{document}
\newtheorem*{theorem}{Theorem}
\newtheorem*{defn}{Definition}
\newtheorem{lemma}{Lemma}
\newtheorem*{prop}{Proposition}
\newtheorem*{cor}{Corollary}

\mbox{ }

\vspace*{-10mm}

\begin{center}
{\Large
Bijections from  Dyck paths to 321-avoiding permutations revisited                         \\ 
}

\vspace{10mm}
DAVID CALLAN  \\
Department of Statistics  \\
\vspace*{-2mm}
University of Wisconsin-Madison  \\
\vspace*{-2mm}
1300 University Ave  \\
\vspace*{-2mm}
Madison, WI \ 53706-1532  \\
{\bf callan@stat.wisc.edu}  \\
\vspace{5mm}

November 16, 2007
\end{center}

\begin{abstract}
There are (at least) three bijections from Dyck paths to 
321-avoiding permutations in the literature, due to Billey-Jockusch-Stanley,
Krattenthaler, and Mansour-Deng-Du. How different are they?
Denoting them $B,K,M$ respectively, we show that $ M = B \circ L = K \circ 
L'$ where $L$ is the classical Kreweras-Lalanne involution on Dyck 
paths and $L'$, also an involution, is a sort of derivative of $L$. Thus 
$K^{-1} \circ B$,
a measure of the difference between $B$ and $K$, is the product of 
involutions $L' \circ L$ and turns out to be a very curious bijection: 
as a permutation on Dyck $n$-paths it  
is an $n$th root of the ``reverse path'' involution. The proof 
of this fact boils down to a geometric argument 
involving pairs of nonintersecting lattice paths.

\end{abstract}

\vspace{5mm}

{\Large \textbf{1 \quad Introduction}  } 
Dyck paths and 321-avoiding permutations are two of the many 
combinatorial manifestations of the Catalan numbers \cite[Ex. 6.19]{ec2}. There are at least three different bijections 
in the literature from Dyck paths to 
321-avoiding permutations, due to  Billey-Jockusch-Stanley 
\cite{bjs} (1993),
Krattenthaler  \cite{kratt2001} (2001) and Mansour-Deng-Du 
\cite{mansour2006} (2006). We denote them $B,K,M$ respectively. (Krattenthaler 
actually gave a bijection to 123-avoiding permutations; $K$ is the 
equivalent bijection to 321-avoiding permutations.) There is 
also a classical involution $L$ on Dyck paths 
dating back to 1970 due to Germain Kreweras \cite{kreweras70} and discussed 
by J.C. Lalanne in 1992-93
\cite{lalanne92, lalanne93}. In this paper we will show that the following 
relationships hold between them:
\begin{equation}
    M = B \circ L = K \circ L'
    \label{eq1}
\end{equation}
where $L'$ is the ``first derivative'' (defined below) of $L$ and, 
like $L$, is an involution. 
We will also see that the bijection $K^{-1} \circ B = L' \circ L$, considered as a 
permutation of Dyck $n$-paths, has order $2n$, a consequence of the fact 
that its $n$th power reverses the path.

The outline of the paper is as follows.
In \S2, we review Dyck path 
terminology and introduce the notion of the derivative $F'$ of a mapping 
$F$ on Dyck paths. Section 3 reviews the ascent-descent code for a 
Dyck path.  Section 4 reviews the left-to-right-maxima and excedance codes for 
a 321-avoiding permutation.  Section 5 describes the involution $L$. 
Section 6 translates $L$ to a simpler setting---pairs of nonintersecting 
lattice paths---and describes $L'$ and $L' \circ L$ in this setting.
Section 7 describes the bijections $B,K$ and $M$.  Section 8 then establishes 
the identities (\ref{eq1}) relating $B,K$ and $M$.
Section 9 uses a geometric argument on nonintersecting path pairs 
to analyze the composition $K^{-1} \circ B =  L' \circ L $.

Astrid Reifegerste \cite{132diagram} has also considered bijections 
involving permutations that avoid a 3-letter pattern and connections between them, 
and some of our observations regarding ``codes'' in \S3 and \S4 can be found in her paper. 
\vspace*{8mm}

{\Large \textbf{2 \quad The derivative of a mapping on Dyck paths}  } 
The set $\DD$ of Dyck paths is the set of lattice paths consisting of 
an equal number of upsteps $u=(1,1)$ and downsteps $d=(1,-1)$ that never dip below 
\gl, the horizontal line connecting its endpoints. The 
\emph{size} or \emph{semilength} of a Dyck path is its number of upsteps. A Dyck $n$-path 
is one of size $n$. An \emph{ascent} is a 
maximal sequence of contiguous upsteps and analogously for a \emph{descent}. 
A \emph{peak} vertex is one preceded by a $u$ and followed by a $d$, and a 
\emph{valley} vertex is defined analogously. 
An \emph{elevated} Dyck path is a nonempty Dyck path whose only return to 
\gl\ occurs at the end. The empty Dyck path is denoted $\ep$. Every nonempty Dyck path decomposes uniquely 
into a concatenation of elevated Dyck paths, called its \emph{components}. For a size-preserving 
bijection $F:\DD \rightarrow \DD$, its \emph{derivative} $F'$ is 
defined by applying $F$ to the ``elevated'' portion of each component 
(and $F'(\ep):=\ep$). 
Schematically,

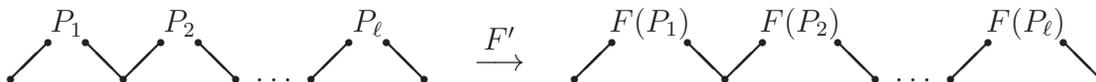
\vspace*{-2mm}

\Einheit=0.5cm
\[
\Pfad(-14,0),3\endPfad
\Pfad(-12,1),43\endPfad
\Pfad(-9,1),4\endPfad
\Pfad(-6,0),3\endPfad
\Pfad(-4,1),4\endPfad
\DuennPunkt(-14,0)
\DuennPunkt(-13,1)
\DuennPunkt(-12,1)
\DuennPunkt(-11,0)
\DuennPunkt(-10,1)
\DuennPunkt(-9,1)
\DuennPunkt(-8,0)
\DuennPunkt(-6,0)
\DuennPunkt(-5,1)
\DuennPunkt(-4,1)
\DuennPunkt(-3,0)
\Label\o{P_{1}}(-12.5,1)
\Label\o{P_{2}}(-9.5,1)
\Label\o{\ldots}(-7,-0.5)
\Label\o{P_{\ell}}(-4.5,1)
\Label\o{F'}(-1,0.6)
\Label\o{\longrightarrow}(-1,-0.1)
\Pfad(1,0),3\endPfad
\Pfad(4,1),43\endPfad
\Pfad(8,1),4\endPfad
\Pfad(11,0),3\endPfad
\Pfad(14,1),4\endPfad
\Label\o{\ldots}(10,-0.5)
\Label\o{F(P_{1})}(3,1)
\Label\o{F(P_{2})}(7,1)
\Label\o{F(P_{\ell})}(13,1)
\DuennPunkt(1,0)
\DuennPunkt(2,1)
\DuennPunkt(4,1)
\DuennPunkt(5,0)
\DuennPunkt(6,1)
\DuennPunkt(8,1)
\DuennPunkt(9,0)
\DuennPunkt(12,1)
\DuennPunkt(14,1)
\DuennPunkt(11,0)
\DuennPunkt(15,0)
\]

Clearly, $F'$ is a
bijection on \DD that preserves not only size but also number of 
components and their sizes, and if $F$ is 
an involution, then so is $F'$.
\vspace*{8mm}

{\Large \textbf{3 \quad The ascent-descent code of a Dyck $\mbf{n}$-path}  }  A 
Dyck path is specified by the lengths of its ascents and descents. For 
example, the path $uuduuuuddduddudd$ has ascent sequence 
$\a=(\a_{i})_{i=1}^{k}=(2,4,1,1)$ and 
descent sequence $\d=(\d_{i})_{i=1}^{k}=(1,3,2,2)$ where $k$ is the 
number of peaks ($ud$s). By definition of Dyck path, each partial sum 
of the ascent lengths 
$\A_{i}:=\sum_{j=1}^{i}\a_{j}$ is $\ge$ the corresponding partial sum 
of the descent lengths $\D_{i}:=\sum_{j=1}^{i}\d_{j}$.
For a Dyck $n$-path, we necessarily have $\A_{k}=\D_{k}=n$ 
and so the path is determined by the pair
$(\A_{i})_{i=1}^{r},\ (\D_{i})_{i=1}^{r}$ where $r:=k-1$ and we call 
this pair the (truncated) \emph{partial-sum ascent-descent code} of the path. 
The preceding example has $n=8,\ k=4,\ r=3$, and 
$\A=(2,6,7),\D=(1,4,6)$.
The precise requirements for a valid partial-sum ascent-descent 
code for a Dyck path of size $n$ are then
\begin{equation}
\begin{array}{ccc}
    0 \le r \le n-1, \\
    1 \le \A_{1} < \A_{2} < \ldots < \A_{r} \le n-1, \\
    1 \le \D_{1} < \D_{2} < \ldots < \D_{r} \le n-1,  \\
   \A_{i}\ge \D_{i}\textrm{ for }1\le i \le r.
\end{array}
\label{cond2}
\end{equation}
Note that the ``pyramid'' path $u^{n}d^{n}$, where exponents denote 
repetition, is the only one with 
$r=0$, and its code consists of two empty sequences.  

\vspace*{8mm}

{\Large \textbf{4 \quad Codes for 321-avoiding permutations}  } A 
permutation $\pi$ on $[n]$ has a left-to-right-maxima decomposition 
as $m_{1}L_{1}m_{2}L_{2}\ldots m_{k}L_{k}$ where 
$m_{1},m_{2},\ldots,m_{k}$ are the left-to-right maxima of $\pi$. For 
example with $n=9$,
\[
\begin{array}{c|cc|c|cc|c|c|c}
    \textrm{{\scriptsize 1}} & \textrm{{\scriptsize 2}} & \textrm{{\scriptsize 3}} & 
    \textrm{{\scriptsize 4}}& \textrm{{\scriptsize 5}} & \textrm{{\scriptsize 
    6}} & 
    \textrm{{\scriptsize 7}} & \textrm{{\scriptsize 8}} & \textrm{{\scriptsize 9}}  \\
    4 & 1 & 3 & 7 & 2 & 5 & 8 & 9 & 6  \\
    m_{1}   & \multicolumn{2}{c|}{L_{1}}  & m_{2} &  
    \multicolumn{2}{c|}{L_{2}} & m_{3} & m_{4} &  
    L_{4}
\end{array}
\]
Here, the left-to-right maxima are $4,7,8,9$ and $L_{3}$ is empty. 
Let's call the left-to-right maxima 
$m=(m_{i})_{i=1}^{k}$ and their positions $p=(p_{i})_{i=1}^{k}$ 
the LRMax skeleton of $\pi$. In the example, $m=(4,7,8,9)$ and 
$p=(1,4,7,8)$. It is easy to see that  
a permutation $\pi$ on $[n]$ is 321-avoiding if and only if the 
concatenated list $L_{1}L_{2}\ldots L_{r}$ is 
increasing. \emph{Thus a $321$-avoiding permutation is determined by its LRMax 
skeleton.}
There are two obvious restrictions on the LRMax 
skeleton: $m_{k}=n$ and $p_{1}=1$. Delete these entries and, to 
make things nice, subtract 1 from each remaining entry in $p$ and 
call the resulting pair, say $(\A_{i})_{i=1}^{r},\ (\D_{i})_{i=1}^{r}$ 
where $r:=k-1$, 
the \emph{LRMax code} of the 321-avoiding permutation on $[n]$.
Note that, for $1\le i \le n-1,\  p_{i+1}\le m_{i}+1$ 
(otherwise the first $m_{i}+1$ entries would all have to be $\le 
m_{i}$, violating the pigeon-hole principle). 
So, since $\A_{i}=m_{i}$ and $\D_{i}=p_{i+1}-1$, the requirements for a valid LRMax code
for a 321-avoiding permutation on $[n]$ are precisely those in 
(\ref{cond2}). This fact is the basis for Krattenthaler's bijection \cite{kratt2001}.

An \emph{excedance location} of a permutation $\pi$ on $[n]$ is an $i\in 
[n-1]$ for which $\pi(i)>i$ (and $\pi(i)$ is then the corresponding 
\emph{excedance value}) and a weak excedance refers to an $i\in 
[n]$ for which $\pi(i)\ge i$. Thus the set of weak excedance 
locations is the disjoint union of the excedance locations and the 
fixed points. Now a 321-avoiding permutation on $[n]$ has the 
following property: if $[n]$ is split into intervals by the fixed 
points $f_{i}$ of $\pi$ so that $[n]$ is the concatenation 
$I_{0}\,f_{1}\,I_{1}\,f_{2}\,\ldots\, f_{q}\,I_{q}$, then $\pi$ preserves each 
interval $I_{i}$. For a 321-avoiding permutation $\pi$ on $[n]$, it follows that 
the left-to-right-maxima coincide with the weak excedance values and that 
the permutation is determined just by its (strict) excedance values $v=(v_{i})$ and 
locations $\ell=(\ell_{i})$. 
In other words, in the LRMax skeleton of a 321-avoiding permutation $\pi$ on $[n]$ 
the fixed points can safely be omitted at the expense of preserving $n$ 
and $\pi^{-1}(n)$ (unless $n$ is a fixed point).
Since each $v_{i}$ is $\ge 2$, let us again subtract 1 to make things 
nice and call the result---$\A:=v-1,\ \D:=\ell$---the \emph{excedance code} 
for $\pi$.
Again, the requirements for a valid excedance code are the same as in 
(\ref{cond2}); this is the basis for the Billey-Jockusch-Stanley 
bijection \cite{bjs}.

\vspace*{8mm}

{\Large \textbf{5 \quad The Lalanne-Kreweras involution on Dyck paths}  } 
We give two descriptions, illustrated with the same example. 

\noindent \textbf{First description \cite{kreweras70,lalanne92,lalanne93}} (graphical): 

\Einheit=0.5cm
\[
\Pfad(-10,4),34343343344433344434\endPfad
\red{\SPfad(-10,4),11111111111111111111\endSPfad}
\SPfad(-5,5),44\endSPfad
\SPfad(-5,5),44\endSPfad
\SPfad(-2,6),44\endSPfad
\SPfad(-3,3),333\endSPfad
\SPfad(0,4),3\endSPfad
\SPfad(3,5),4\endSPfad
\SPfad(4,6),44\endSPfad
\SPfad(4,4),33\endSPfad
\SPfad(6,4),3\endSPfad
\DuennPunkt(-10,4)
\DuennPunkt(-9,5)
\DuennPunkt(-8,4)
\DuennPunkt(-7,5)
\DuennPunkt(-6,4)
\DuennPunkt(-5,5)
\DuennPunkt(-4,6)
\DuennPunkt(-3,5)
\DuennPunkt(-2,6)
\DuennPunkt(-1,7)
\DuennPunkt(0,6)
\DuennPunkt(1,5)
\DuennPunkt(2,4)
\DuennPunkt(3,5)
\DuennPunkt(4,6)
\DuennPunkt(5,7)
\DuennPunkt(6,6)
\DuennPunkt(7,5)
\DuennPunkt(8,4)
\DuennPunkt(9,5)
\DuennPunkt(10,4)
\Label\o{\textrm{\footnotesize 1}}(-5.5,4.5)
\Label\o{\textrm{\footnotesize 2}}(-2.5,5.5)
\Label\o{\textrm{\footnotesize 3}}(2.5,4.5)
\Label\o{\textrm{\footnotesize 4}}(3.5,5.5)
\blue{
\Label\o{\textrm{\footnotesize 1}}(-3,1.8)
\Label\o{\textrm{\footnotesize 2}}(0,2.8)
\Label\o{\textrm{\footnotesize 3}}(4,2.8)
\Label\o{\textrm{\footnotesize 4}}(6,2.8)
\DickPunkt(-3,3)
\DickPunkt(0,4)
\DickPunkt(4,4)
\DickPunkt(6,4)
\DuennPunkt(-10,4)
\DuennPunkt(-9,3)
\DuennPunkt(-8,2)
\DuennPunkt(-7,1)
\DuennPunkt(-6,0)
\DuennPunkt(-5,1)
\DuennPunkt(-4,2)
\DuennPunkt(-3,3)
\DuennPunkt(-2,2)
\DuennPunkt(-1,3)
\DuennPunkt(0,4)
\DuennPunkt(1,3)
\DuennPunkt(2,2)
\DuennPunkt(3,3)
\DuennPunkt(4,4)
\DuennPunkt(5,3)
\DuennPunkt(6,4)
\DuennPunkt(7,3)
\DuennPunkt(8,2)
\DuennPunkt(9,3)
\DuennPunkt(10,4)
\Pfad(-10,4),44443334334433434433\endPfad }
\]
Draw a southeast line from the midpoint of each $uu$ and a 
southwest line from the midpoint of each $dd$. There will be the same number 
of each. Mark the point of intersection of the $i$th  southeast and 
the $i$th southwest 
line for each $i$. Then form the unique (inverted) Dyck path with (inverted)
valleys at the marked points, as shown in blue (below \gl) above.

\noindent \textbf{Second description} (algorithmic): 
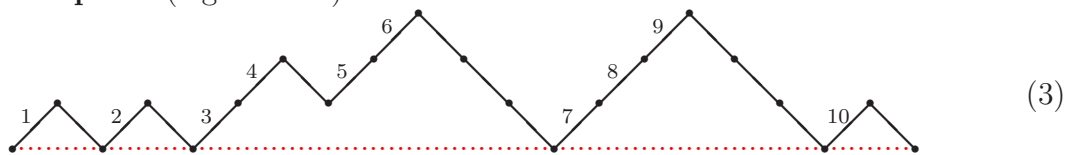
\Einheit=0.6cm
\begin{equation}
\Pfad(-10,-1),34343343344433344434\endPfad
\red{\SPfad(-10,-1),11111111111111111111\endSPfad}
\Label\o{\textrm{\scriptsize 1}}(-9.7,-0.7)
\Label\o{\textrm{\scriptsize 2}}(-7.7,-0.7)
\Label\o{\textrm{\scriptsize 3}}(-5.7,-0.7)
\Label\o{\textrm{\scriptsize 4}}(-4.7,0.3)
\Label\o{\textrm{\scriptsize 5}}(-2.7,0.3)
\Label\o{\textrm{\scriptsize 6}}(-1.7,1.3)
\Label\o{\textrm{\scriptsize 7}}(2.3,-0.7)
\Label\o{\textrm{\scriptsize 8}}(3.3,0.3)
\Label\o{\textrm{\scriptsize 9}}(4.3,1.3)
\Label\o{\textrm{\scriptsize 10}}(8.3,-0.7)
\DuennPunkt(-10,-1)
\DuennPunkt(-9,0)
\DuennPunkt(-8,-1)
\DuennPunkt(-7,0)
\DuennPunkt(-6,-1)
\DuennPunkt(-5,0)
\DuennPunkt(-4,1)
\DuennPunkt(-3,0)
\DuennPunkt(-2,1)
\DuennPunkt(-1,2)
\DuennPunkt(0,1)
\DuennPunkt(1,0)
\DuennPunkt(2,-1)
\DuennPunkt(3,0)
\DuennPunkt(4,1)
\DuennPunkt(5,2)
\DuennPunkt(6,1)
\DuennPunkt(7,0)
\DuennPunkt(8,-1)
\DuennPunkt(9,0)
\DuennPunkt(10,-1)
\label{dyck1}
\end{equation}
\vspace*{1mm}

\noindent Label the upsteps left to right. Record the label on the first 
$u$ of each occurrence of $uu$. 
The example gives $(3,5,7,8)$. Call this vector \D. Do likewise for the downsteps. The example 
gives $(4,5,7,8)$. Call this vector \A. Then form the Dyck path whose partial-sum 
ascent-descent code is $(\A,\D)$. (The reader may check that $(\A,\D)$ 
satisfy the defining conditions (\ref{cond2}) with $n$ the size of the 
path). The example gives ascent lengths 
4,1,2,1,2 and descent lengths 3,2,2,1,2.

In the next section, following Emeric Deutsch \cite{em99}, we use a 
suitable bijection to identify Dyck paths with another manifestation 
of the Catalan numbers, nonintersecting path pairs (parallelogram 
polyominoes). In this setting $L$ has perhaps its simplest possible 
description: flip the path pair in a $45^{\circ}$ line. Also, the 
``reverse path'' involution $R$ on Dyck paths translates to ``rotate path pair $180^{\circ}$''.

\vspace*{8mm}

{\Large \textbf{6 \quad $\mbf{L,\ L'},\ \mbf{L}\circ \mbf{L'}$, and $\mbf{R}$ on Path Pairs}  } 

A \emph{$($nonintersecting$)$ path pair} is an ordered pair $(P_{1},P_{2})$ of 
paths of unit steps north, $N=(1,0)$, and east, $E=(0,1)$, that intersect only at the 
initial and terminal points and such that $P_{1}$ (the upper path) lies above 
$P_{2}$. The size of a path pair is the number of steps in each path, necessarily 
the same. The region enclosed by a path pair is known as a 
parallelogram polyomino.

There is a well known bijection 
\cite[p.\,182]{delest84}\cite[Ex. 6.19($\ell$)]{ec2} which we will use to 
identify Dyck paths of size $n$ with 
path pairs of size $n+1$. An equivalent bijection (up to reversing Dyck 
paths and rotating path pairs) has been given by Sulanke \cite[p.\,295]{sulanke93}.
Here is the bijection (with a slightly simplified description).

Given a Dyck path, first elevate it, that is, prepend $u$ and append $d$.

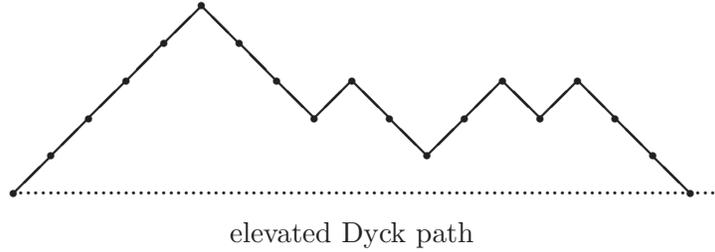
\Einheit=0.5cm
\[
\Pfad(-9,0),333334443443343444\endPfad
\SPfad(-9,0),1111111111111111111\endSPfad
\DuennPunkt(-9,0)
\DuennPunkt(-8,1)
\DuennPunkt(-7,2)
\DuennPunkt(-6,3)
\DuennPunkt(-5,4)
\DuennPunkt(-4,5)
\DuennPunkt(-3,4)
\DuennPunkt(-2,3)
\DuennPunkt(-1,2)
\DuennPunkt(0,3)
\DuennPunkt(1,2)
\DuennPunkt(2,1)
\DuennPunkt(3,2)
\DuennPunkt(4,3)
\DuennPunkt(5,2)
\DuennPunkt(6,3)
\DuennPunkt(7,2)
\DuennPunkt(8,1)
\DuennPunkt(9,0)
\Label\u{ \textrm{{\small elevated Dyck path}}}(0,-0.5)
\]

\vspace*{4mm}

Then extract the elevated path's ascents as $N$ steps except that the last step in each 
ascent is rendered as an $E$ step, and concatenate:
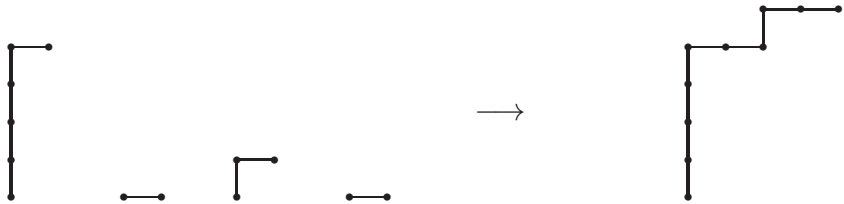
\Einheit=0.5cm
\[
\Pfad(-11,0),22221\endPfad
\Pfad(-8,0),1\endPfad
\Pfad(-5,0),21\endPfad
\Pfad(-2,0),1\endPfad
\DuennPunkt(-11,0)
\DuennPunkt(-11,1)
\DuennPunkt(-11,2)
\DuennPunkt(-11,3)
\DuennPunkt(-11,4)
\DuennPunkt(-10,4)
\DuennPunkt(-8,0)
\DuennPunkt(-7,0)
\DuennPunkt(-5,0)
\DuennPunkt(-5,1)
\DuennPunkt(-4,1)
\DuennPunkt(-2,0)
\DuennPunkt(-1,0)
\Pfad(7,0),222211211\endPfad
\DuennPunkt(7,0)
\DuennPunkt(7,1)
\DuennPunkt(7,2)
\DuennPunkt(7,3)
\DuennPunkt(7,4)
\DuennPunkt(8,4)
\DuennPunkt(9,4)
\DuennPunkt(9,5)
\DuennPunkt(10,5)
\DuennPunkt(11,5)
\Label\u{ \longrightarrow}(2,2.8)
\]

\vspace*{2mm}

This is the upper path.

Do likewise for the descents to get a path, $X$ say, and then 
\emph{transfer the last step, necessarily an $E$, to the start}. This 
gives the lower path and the resulting path pair is
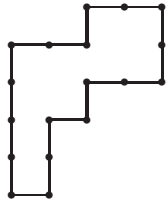
\Einheit=0.5cm
\[
\Pfad(-2,0),222211211\endPfad
\Pfad(-2,0),122121122\endPfad
\DuennPunkt(-2,0)
\DuennPunkt(-2,1)
\DuennPunkt(-2,2)
\DuennPunkt(-2,3)
\DuennPunkt(-2,4)
\DuennPunkt(-1,4)
\DuennPunkt(0,4)
\DuennPunkt(0,5)
\DuennPunkt(1,5)
\DuennPunkt(2,5)
\DuennPunkt(-1,0)
\DuennPunkt(-1,1)
\DuennPunkt(-1,2)
\DuennPunkt(-0,2)
\DuennPunkt(0,3)
\DuennPunkt(1,3)
\DuennPunkt(2,3)
\DuennPunkt(2,4)
\]

If we let $\a_{i}$ denote the length of the $i$th ascent in the 
elevated Dyck path and 
$\d_{i}$ the length of the $i$th descent for $1\le i \le k$, where 
$k=\#\,$ peaks (= \#\,ascents = \#\,descents), then, since the path 
is elevated, 
$\sum_{i=1}^{j}\a_{i}>\sum_{i=1}^{j}\d_{i}$ for $j=1,2,\ldots,k-1$, and hence the $i$th 
$E$ step in the upper path lies strictly above the $i$th $E$ step in 
the path $X$ for $i< k$. This ensures that the resulting path pair is nonintersecting 
and the mapping is clearly invertible. Let us call this bijection $\phi$.

Using $\phi$ to identify Dyck paths and nonintersecting path 
pairs, the Kreweras-Lalanne involution $L$ simplifies to ``flip path 
pair in a $45^{\circ}$ line''. Again using the Dyck path (\ref{dyck1}) 
from \S 5 to illustrate,

\Einheit=0.5cm
\[
\Pfad(-8,0),21121212211\endPfad
\Pfad(-8,0),11112212212\endPfad
\Pfad(2,0),22221121121\endPfad
\Pfad(2,0),12212121122\endPfad
\Label\o{\longrightarrow}(0,2)
\Label\o{\longrightarrow}(-9.5,2)
\Label\o{\longrightarrow}(9.5,2)
\DuennPunkt(-8,0)
\DuennPunkt(-7,0)
\DuennPunkt(-8,1)
\DuennPunkt(-7,1)
\DuennPunkt(-6,0)
\DuennPunkt(-5,0)
\DuennPunkt(-4,0)
\DuennPunkt(-3,2)
\DuennPunkt(-2,4)
\DuennPunkt(-6,1)
\DuennPunkt(-5,2)
\DuennPunkt(-4,3)
\DuennPunkt(-3,3)
\DuennPunkt(-2,5)
\DuennPunkt(-6,2)
\DuennPunkt(-5,3)
\DuennPunkt(-4,4)
\DuennPunkt(-4,5)
\DuennPunkt(-3,4)
\DuennPunkt(-3,5)
\DuennPunkt(2,4)
\DuennPunkt(2,3)
\DuennPunkt(2,2)
\DuennPunkt(2,1)
\DuennPunkt(2,0)
\DuennPunkt(3,4)
\DuennPunkt(3,2)
\DuennPunkt(3,1)
\DuennPunkt(3,0)
\DuennPunkt(4,5)
\DuennPunkt(4,4)
\DuennPunkt(4,3)
\DuennPunkt(4,2)
\DuennPunkt(5,5)
\DuennPunkt(6,6)
\DuennPunkt(7,6)
\DuennPunkt(5,4)
\DuennPunkt(6,5)
\DuennPunkt(7,5)
\DuennPunkt(5,3)
\DuennPunkt(6,4)
\DuennPunkt(7,4)
\Label\o{2\ 1\ 2\ 2\ 3\ 1}(-13,1)
\Label\o{1\ 1\ 1\ 3\ 3\ 2}(-13,0)
\Label\u{\textrm{{\footnotesize elevated}}}(-13,-1)
\Label\u{\textrm{{\footnotesize Dyck path}}}(-13,-1.8)
\Label\u{\textrm{{\footnotesize ascent/descent}}}(-13,-2.6)
\Label\u{\textrm{{\footnotesize lengths}}}(-13,-3.4)
\Label\u{\textrm{{\footnotesize path pair}}}(-5,-1)
\Label\u{\textrm{{\footnotesize path pair}}}(5,-1)
\Label\u{\textrm{{\footnotesize flip}}}(0,2)
\Label\o{5\ 1\ 2\ 1\ 2}(13,1)
\Label\o{3\ 2\ 2\ 1\ 3}(13,0)
\Label\u{\textrm{{\footnotesize elevated}}}(13,-1)
\Label\u{\textrm{{\footnotesize Dyck path}}}(13,-1.8)
\Label\u{\textrm{{\footnotesize ascent/descent}}}(13,-2.6)
\Label\u{\textrm{{\footnotesize lengths}}}(13,-3.4)
\]

\vspace*{2mm}

To see the effect of $L'$ on a path pair $P=\phi(D)$ where $D$ is a 
Dyck path, we need to 
identify within $P$ the interior of each component of $D$ (in blue 
below), and this is easy 
to do. 

\Einheit=0.5cm
\[
\Pfad(-13,0),33344344343333443444333444\endPfad
\SPfad(-13,0),11111111111111111111111111\endSPfad 
\blue{\Pfad(-12,1),334434\endPfad
\Pfad(-2,1),33344344\endPfad
\Pfad(8,1),3344\endPfad
\DuennPunkt(-13,0)
\DuennPunkt(-12,1)
\DuennPunkt(-11,2)
\DuennPunkt(-10,3)
\DuennPunkt(-9,2)
\DuennPunkt(-8,1)
\DuennPunkt(-7,2)
\DuennPunkt(-6,1)
\DuennPunkt(-5,0)
\NormalPunkt(-4,1)
\DuennPunkt(-3,0)
\DuennPunkt(-2,1)
\DuennPunkt(-1,2)
\DuennPunkt(0,3)
\DuennPunkt(1,4)
\DuennPunkt(2,3)
\DuennPunkt(3,2)
\DuennPunkt(4,3)
\DuennPunkt(5,2)
\DuennPunkt(6,1)
\DuennPunkt(7,0)
\DuennPunkt(8,1)
\DuennPunkt(9,2)
\DuennPunkt(10,3)
\DuennPunkt(11,2)
\DuennPunkt(12,1)
\DuennPunkt(13,0)}
\blue{\SPfad(-12,1),111111\endSPfad
\SPfad(-2,1),11111111\endSPfad
\SPfad(8,1),1111\endSPfad}
\DuennPunkt(-13,0)
\DuennPunkt(-5,0)
\DuennPunkt(-3,0)
\DuennPunkt(7,0)
\DuennPunkt(13,0)
\Label\o{{\textrm{\small Dyck path $D$}}}(0,-1.5)
\]

\vspace*{3mm}

\noindent The components 
of $D$ are determined by the points of contact of $D$ with \gl. These 
points, including the initial and terminal points, correspond to 
unit vertical segments (in red in the figure below right) joining a vertex of 
the upper path to a vertex of the lower path. Furthermore, each hill ($ud$ pair at \gl) in 
$D$ corresponds to a pair of steps in $P$ that 
form horizontal sides of a unit square. Keep in mind that $\phi$ sends 
a Dyck $n$-path  to a path pair of size $n+1$ \emph{except when $n=0$}: the 
empty Dyck path corresponds to the empty path pair. So we can expect 
a hill in $D$---a component with empty interior---to exhibit singular 
behavior under $L'$. Indeed, the path pair corresponding to the 
interior of each component of $D$ can be seen in $P=\phi(D)$ as in the 
illustration, where numerals label steps in each upper path and letters 
in each lower path, and hills in $D$ show up in $P$ as unlabeled unit squares 
bounded above by $P_{1}$ and below by $P_{2}$.

\vspace{1mm}

\Einheit=0.6cm
\[
\Pfad(-11,0),2211\endPfad
\Pfad(-11,0),1212\endPfad
\Pfad(-4,0),22211\endPfad
\Pfad(-4,0),12122\endPfad
\Pfad(0,0),221\endPfad
\Pfad(0,0),122\endPfad
\Pfad(5,0),22211122211221\endPfad
\Pfad(5,0),12121121221222\endPfad
\red{
\Pfad(5,0),2\endPfad
\Pfad(7,2),2\endPfad
\Pfad(8,2),2\endPfad
\Pfad(10,5),2\endPfad
\Pfad(11,7),2\endPfad }
\Label\o{\textrm{{\footnotesize The ascent lengths \a and descent 
lengths \d of}}}(-6,9.5)
\Label\o{\textrm{{\footnotesize the components of $D$ are just 
what is needed}}}(-6,8.7)
\Label\o{\textrm{{\footnotesize to construct the path pairs for their 
interiors}}}(-6,7.9)
\Label\o{\emptyset}(-6.5,0)
\Label\o{\d:}(-12.5,4.5)
\Label\o{\a:}(-12.5,5.3)
\Label\o{{\textrm{\footnotesize 2}}}(-9.5,4.5)
\Label\o{{\textrm{\footnotesize 2}}}(-10.5,4.5)
\Label\o{{\textrm{\footnotesize 1}}}(-9.5,5.3)
\Label\o{{\textrm{\footnotesize 3}}}(-10.5,5.3)
\Label\o{{\textrm{\footnotesize 1}}}(-6.5,4.5)
\Label\o{{\textrm{\footnotesize 1}}}(-6.5,5.3)
\Label\o{{\textrm{\footnotesize 3}}}(-2.5,4.5)
\Label\o{{\textrm{\footnotesize 2}}}(-3.5,4.5)
\Label\o{{\textrm{\footnotesize 1}}}(-2.5,5.3)
\Label\o{{\textrm{\footnotesize 4}}}(-3.5,5.3)
\Label\o{{\textrm{\footnotesize 3}}}(0.5,4.5)
\Label\o{{\textrm{\footnotesize 3}}}(0.5,5.3)
\Label\o{{\textrm{\small path pair $P=\phi(D)$}}}(8,9.5)
\Label\o{{\textrm{\scriptsize 1}}}(-11.3,0.1)
\Label\o{{\textrm{\scriptsize 2}}}(-11.3,1.1)
\Label\o{{\textrm{\scriptsize 3}}}(-10.5,2)
\Label\o{{\textrm{\scriptsize 4}}}(-9.5,2)
\Label\u{{\textrm{\scriptsize $a$}}}(-10.5,0.2)
\Label\o{{\textrm{\scriptsize $b$}}}(-10.2,0.1)
\Label\u{{\textrm{\scriptsize $c$}}}(-9.5,1.2)
\Label\o{{\textrm{\scriptsize $d$}}}(-8.7,1.1)
\Label\o{{\textrm{\scriptsize 5}}}(-4.3,0.1)
\Label\o{{\textrm{\scriptsize 6}}}(-4.3,1.1)
\Label\o{{\textrm{\scriptsize 7}}}(-4.3,2.1)
\Label\o{{\textrm{\scriptsize 8}}}(-3.5,3)
\Label\o{{\textrm{\scriptsize 9}}}(-2.5,3)
\Label\u{{\textrm{\scriptsize $e$}}}(-3.5,0.2)
\Label\o{{\textrm{\scriptsize $f$}}}(-3.2,0.1)
\Label\u{{\textrm{\scriptsize $g$}}}(-2.5,1.2)
\Label\o{{\textrm{\scriptsize $h$}}}(-1.7,1.1)
\Label\o{{\textrm{\scriptsize $i$}}}(-1.7,2.1)
\Label\o{{\textrm{\scriptsize 10}}}(-0.4,0.1)
\Label\o{{\textrm{\scriptsize 11}}}(-0.4,1.1)
\Label\o{{\textrm{\scriptsize 12}}}(0.5,2)
\Label\u{{\textrm{\scriptsize $j$}}}(0.5,0.2)
\Label\o{{\textrm{\scriptsize $k$}}}(1.3,0.1)
\Label\o{{\textrm{\scriptsize $\ell$}}}(1.3,1.1)
\DuennPunkt(-11,0)
\DuennPunkt(-11,1)
\DuennPunkt(-11,2)
\DuennPunkt(-10,0)
\DuennPunkt(-10,1)
\DuennPunkt(-10,2)
\DuennPunkt(-9,1)
\DuennPunkt(-9,2)
\DuennPunkt(-4,0)
\DuennPunkt(-4,1)
\DuennPunkt(-4,2)
\DuennPunkt(-4,3)
\DuennPunkt(-3,0)
\DuennPunkt(-3,1)
\DuennPunkt(-3,3)
\DuennPunkt(-2,1)
\DuennPunkt(-2,2)
\DuennPunkt(-2,3)
\DuennPunkt(0,0)
\DuennPunkt(0,1)
\DuennPunkt(0,2)
\DuennPunkt(1,0)
\DuennPunkt(1,1)
\DuennPunkt(1,2)
\DuennPunkt(5,0)
\DuennPunkt(5,1)
\DuennPunkt(5,2)
\DuennPunkt(5,3)
\DuennPunkt(6,0)
\DuennPunkt(6,1)
\DuennPunkt(6,3)
\DuennPunkt(7,1)
\DuennPunkt(7,2)
\DuennPunkt(7,3)
\DuennPunkt(8,2)
\DuennPunkt(8,3)
\DuennPunkt(8,4)
\DuennPunkt(8,5)
\DuennPunkt(8,6)
\DuennPunkt(9,2)
\DuennPunkt(9,3)
\DuennPunkt(9,6)
\DuennPunkt(10,3)
\DuennPunkt(10,4)
\DuennPunkt(10,5)
\DuennPunkt(10,6)
\DuennPunkt(10,7)
\DuennPunkt(10,8)
\DuennPunkt(11,5)
\DuennPunkt(11,6)
\DuennPunkt(11,7)
\DuennPunkt(11,8)
\Label\o{{\textrm{\scriptsize 1}}}(4.7,1.1)
\Label\o{{\textrm{\scriptsize 2}}}(4.7,2.1)
\Label\o{{\textrm{\scriptsize 3}}}(5.5,3)
\Label\o{{\textrm{\scriptsize 4}}}(6.5,3)
\Label\o{{\textrm{\scriptsize 8}}}(8.5,6)
\Label\o{{\textrm{\scriptsize 9}}}(9.5,6)
\Label\o{{\textrm{\scriptsize 12}}}(10.5,8)
\Label\o{{\textrm{\scriptsize 5}}}(7.7,3.1)
\Label\o{{\textrm{\scriptsize 6}}}(7.7,4.1)
\Label\o{{\textrm{\scriptsize 7}}}(7.7,5.1)
\Label\o{{\textrm{\scriptsize 10}}}(10.3,6.1)
\Label\o{{\textrm{\scriptsize 11}}}(10.3,7.1)
\Label\u{{\textrm{\scriptsize $a$}}}(5.5,0.2)
\Label\o{{\textrm{\scriptsize $b$}}}(5.8,0.1)
\Label\u{{\textrm{\scriptsize $c$}}}(6.5,1.2)
\Label\o{{\textrm{\scriptsize $d$}}}(6.8,1.1)
\Label\u{{\textrm{\scriptsize $e$}}}(8.5,2.2)
\Label\o{{\textrm{\scriptsize $f$}}}(8.8,2.1)
\Label\u{{\textrm{\scriptsize $g$}}}(9.5,3.2)
\Label\o{{\textrm{\scriptsize $h$}}}(9.8,3.1)
\Label\o{{\textrm{\scriptsize $i$}}}(9.8,4.1)
\Label\u{{\textrm{\scriptsize $j$}}}(10.5,5.2)
\Label\o{{\textrm{\scriptsize $k$}}}(11.3,5.1)
\Label\o{{\textrm{\scriptsize $\ell$}}}(11.3,6.1)
\]

\vspace*{2mm}

 Since $L$ flips a path pair, 
the effect of $L'$ is to flip in a $45^{\circ}$ line the pairs of 
labeled segments separated 
by red lines while preserving the red lines and unlabeled unit squares. 
The example yields 

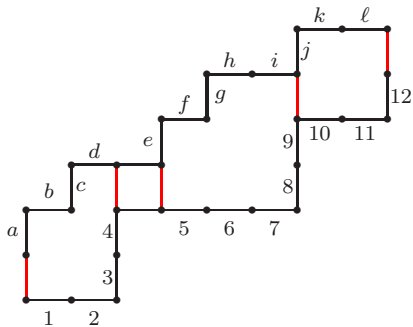
\Einheit=0.6cm
\[
\Pfad(-4,0),22121121211211\endPfad
\Pfad(-4,0),11221111221122\endPfad
\red{
\Pfad(-4,0),2\endPfad
\Pfad(-2,2),2\endPfad
\Pfad(-1,2),2\endPfad
\Pfad(2,4),2\endPfad
\Pfad(4,5),2\endPfad }
\DuennPunkt(-4,0)
\DuennPunkt(-4,1)
\DuennPunkt(-4,2)
\DuennPunkt(-3,0)
\DuennPunkt(-3,2)
\DuennPunkt(-3,3)
\DuennPunkt(-2,0)
\DuennPunkt(-2,1)
\DuennPunkt(-2,2)
\DuennPunkt(-2,3)
\DuennPunkt(-1,2)
\DuennPunkt(-1,3)
\DuennPunkt(-1,4)
\DuennPunkt(0,2)
\DuennPunkt(0,4)
\DuennPunkt(0,5)
\DuennPunkt(1,2)
\DuennPunkt(1,5)
\DuennPunkt(2,2)
\DuennPunkt(2,3)
\DuennPunkt(2,4)
\DuennPunkt(2,5)
\DuennPunkt(2,6)
\DuennPunkt(3,4)
\DuennPunkt(3,6)
\DuennPunkt(4,4)
\DuennPunkt(4,5)
\DuennPunkt(4,6)
\Label\o{{\textrm{\scriptsize $a$}}}(-4.3,1.1)
\Label\o{{\textrm{\scriptsize $b$}}}(-3.5,2.0)
\Label\o{{\textrm{\scriptsize $c$}}}(-2.8,2.1)
\Label\o{{\textrm{\scriptsize $d$}}}(-2.5,2.9)
\Label\o{{\textrm{\scriptsize $e$}}}(-1.3,3.1)
\Label\o{{\textrm{\scriptsize $f$}}}(-0.5,3.9)
\Label\o{{\textrm{\scriptsize $g$}}}(0.3,4.1)
\Label\o{{\textrm{\scriptsize $h$}}}(0.5,4.9)
\Label\o{{\textrm{\scriptsize $i$}}}(1.5,4.9)
\Label\o{{\textrm{\scriptsize $j$}}}(2.2,5.1)
\Label\o{{\textrm{\scriptsize $k$}}}(2.5,5.9)
\Label\o{{\textrm{\scriptsize $\ell$}}}(3.5,5.9)
\Label\u{{\textrm{\scriptsize $1$}}}(-3.5,0.1)
\Label\u{{\textrm{\scriptsize $2$}}}(-2.5,0.1)
\Label\o{{\textrm{\scriptsize $3$}}}(-2.2,0.1)
\Label\o{{\textrm{\scriptsize $4$}}}(-2.2,1.1)
\Label\u{{\textrm{\scriptsize $5$}}}(-0.5,2.1)
\Label\u{{\textrm{\scriptsize $6$}}}(0.5,2.1)
\Label\u{{\textrm{\scriptsize $7$}}}(1.5,2.1)
\Label\o{{\textrm{\scriptsize $8$}}}(1.8,2.1)
\Label\o{{\textrm{\scriptsize $9$}}}(1.8,3.1)
\Label\u{{\textrm{\scriptsize $10$}}}(2.5,4.2)
\Label\u{{\textrm{\scriptsize $11$}}}(3.5,4.2)
\Label\o{{\textrm{\scriptsize $12$}}}(4.3,4.1)
\Label\o{{\textrm{\small $L'(P)$}}}(0,-2)
\]

\vspace*{2mm}

\noindent where labels and interior red lines are included for clarity.

Now, since $L$ flips the entire path pair and $L'$ flips a good deal 
of it, the composition $L'\circ L$ merely tweaks it: 
given a path pair $P$, to obtain $L'\circ L(P)$

\noindent (i) identify the last (northeasternmost) step in the upper 
path and the southwesternmost step in the lower path. Both these steps 
are necessarily flat.

\noindent (ii) identify each pair of \emph{vertical} steps that form opposite 
sides of a unit square. 

Then change the two steps in (i) from flat to vertical and all steps 
in (ii) (if any) from vertical to flat. Two examples are shown below 
(unchanged steps in color, labels for clarity). Note that step (ii) ensures the resulting path pair is nonintersecting.

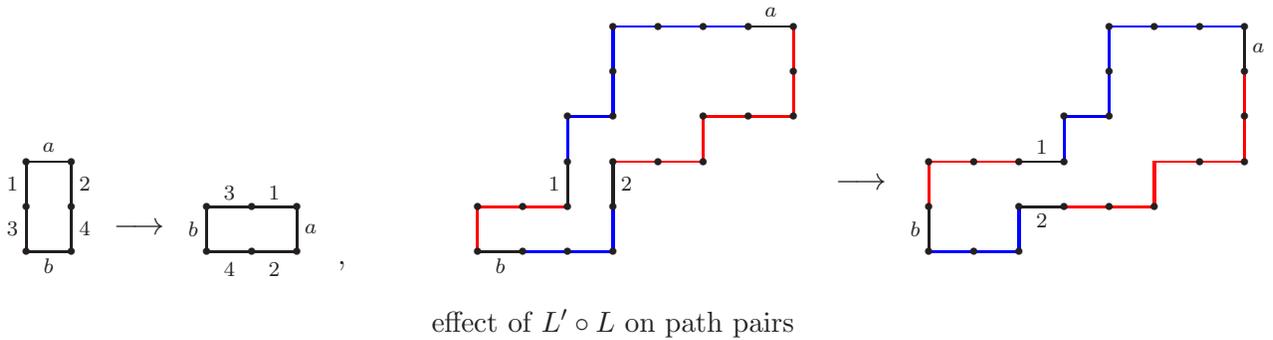
\vspace{-2mm}

\Einheit=0.6cm
\[
\Pfad(-13,0),221\endPfad
\Pfad(-13,0),122\endPfad
\Pfad(-9,0),211\endPfad
\Pfad(-9,0),112\endPfad
\Pfad(-3,0),211221221111\endPfad
\Pfad(-3,0),111221121122\endPfad
\Pfad(7,0),221112122111\endPfad
\Pfad(7,0),112111211222\endPfad
\red{
\Pfad(-3,0),211\endPfad
\Pfad(0,2),1121122\endPfad
\Pfad(7,1),211\endPfad
\Pfad(10,1),1121122\endPfad
 }
\blue{\Pfad(-1,2),2122111\endPfad
 \Pfad(10,2),2122111\endPfad
 \Pfad(-2,0),112\endPfad
 \Pfad(7,0),112\endPfad }
\DuennPunkt(-13,0)
\DuennPunkt(-13,1)
\DuennPunkt(-13,2)
\DuennPunkt(-12,0)
\DuennPunkt(-12,1)
\DuennPunkt(-12,2)
\DuennPunkt(-9,0)
\DuennPunkt(-9,1)
\DuennPunkt(-8,0)
\DuennPunkt(-8,1)
\DuennPunkt(-7,0)
\DuennPunkt(-7,1) 
\DuennPunkt(-3,0)
\DuennPunkt(-3,1)
\DuennPunkt(-2,0)
\DuennPunkt(-2,1)
\DuennPunkt(-1,0)
\DuennPunkt(-1,1)
\DuennPunkt(-1,2)
\DuennPunkt(-1,3)
\DuennPunkt(0,0)
\DuennPunkt(0,1)
\DuennPunkt(0,2)
\DuennPunkt(0,3)
\DuennPunkt(0,4)
\DuennPunkt(0,5)
\DuennPunkt(1,2)
\DuennPunkt(1,5)
\DuennPunkt(2,2)
\DuennPunkt(2,3)
\DuennPunkt(2,5)
\DuennPunkt(3,3)
\DuennPunkt(3,5)
\DuennPunkt(4,3)
\DuennPunkt(4,4)
\DuennPunkt(4,5)
\DuennPunkt(7,1)
\DuennPunkt(7,2)
\DuennPunkt(7,0)
\DuennPunkt(8,2)
\DuennPunkt(8,0)
\DuennPunkt(9,2)
\DuennPunkt(9,1)
\DuennPunkt(9,0)
\DuennPunkt(10,3)
\DuennPunkt(10,2)
\DuennPunkt(10,1)
\DuennPunkt(11,5)
\DuennPunkt(11,4)
\DuennPunkt(11,3)
\DuennPunkt(11,1)
\DuennPunkt(12,1)
\DuennPunkt(12,5)
\DuennPunkt(13,2)
\DuennPunkt(13,5)
\DuennPunkt(14,2)
\DuennPunkt(14,3)
\DuennPunkt(14,4)
\DuennPunkt(14,5)
\Label\o{{\textrm{\scriptsize 3}}}(-13.3,0.1)
\Label\o{{\textrm{\scriptsize 1}}}(-13.3,1.1)
\Label\o{{\textrm{\scriptsize 2}}}(-11.7,1.1)
\Label\o{{\textrm{\scriptsize 4}}}(-11.7,0.1)
\Label\u{{\textrm{\scriptsize $b$}}}(-12.5,0.2)
\Label\o{{\textrm{\scriptsize $a$}}}(-12.5,1.9)
\Label\o{{\textrm{\scriptsize $b$}}}(-9.3,0.1)
\Label\o{{\textrm{\scriptsize $a$}}}(-6.7,0.1)
\Label\o{{\textrm{\scriptsize 3}}}(-8.5,0.9)
\Label\o{{\textrm{\scriptsize 1}}}(-7.5,0.9)
\Label\u{{\textrm{\scriptsize 2}}}(-7.5,0.1)
\Label\u{{\textrm{\scriptsize 4}}}(-8.5,0.1)
\Label\o{\longrightarrow}(-10.5,0.1)
\Label\u{{\textrm{\scriptsize $b$}}}(-2.5,0.2)
\Label\o{{\textrm{\scriptsize $a$}}}(3.5,4.9)
\Label\o{{\textrm{\scriptsize 1}}}(-1.3,1.1)
\Label\o{{\textrm{\scriptsize 2}}}(0.3,1.1)
\Label\u{{\textrm{\scriptsize $2$}}}(9.5,1.2)
\Label\o{{\textrm{\scriptsize $a$}}}(14.3,4.1)
\Label\o{{\textrm{\scriptsize 1}}}(9.5,1.9)
\Label\o{{\textrm{\scriptsize $b$}}}(6.7,0.1)
\Label\u{,}(-6,0.2)
\Label\o{\longrightarrow}(5.5,1.1)
\Label\o{{\textrm{\small effect of $L'\circ L$ on path pairs}}}(0,-2.0)
\]
\vspace{1mm}

Finally, it is clear that reversing a Dyck path, which interchanges 
the roles of upsteps and downsteps, corresponds under $\phi$ to 
rotating a path pair $180^{\circ}$.

\vspace*{8mm}

{\Large \textbf{7 \quad The bijections $\mbf{B,K}$  and $\mbf{M}$}  } 
The Billey-Jockusch-Stanley bijection $B$ from Dyck paths to 321-avoiding 
permutations can now be simply described: form the partial-sum ascent-descent 
code $(\A,\D)$ of the Dyck path and then use it as the excedance code of a 321-avoiding 
permutation. For example, the Dyck path (\ref{dyck1}) of \S 5
has size $n=10$, ascent lengths $\a=(1,1,2,2,3,1)$ and descent lengths 
$\d=(1,1,1,3,3,1)$ so that
$\A=(1,2,4,6,9),\ \D=(1,2,3,6,9)$. With $(\A,\D)$ as excedance 
code, $\A+1$ gives the excedance values and \D the excedance 
locations. We thus immediately have the following partial permutation
\[
\begin{array}{cccccccccc}
 \textrm{{\scriptsize 1}} & \textrm{{\scriptsize 2}} & \textrm{{\scriptsize 3}} & 
    \textrm{{\scriptsize 4}}& \textrm{{\scriptsize 5}} & \textrm{{\scriptsize 
    6}} & 
    \textrm{{\scriptsize 7}} & \textrm{{\scriptsize 8}} 
    & \textrm{{\scriptsize 9}} & \textrm{{\scriptsize 10}}   \\[-0.8ex]
    2 & 3 & 5 &   &  & 7 &  &  & 10 & 
\end{array}
\]
and filling in the missing entries in increasing order gives the 
image permutation: \linebreak
$2\ 3\ 5\ 1\ 4\ 7\ 6\ 8\ 10\ 9$.

The Krattenthaler bijection $K$ uses the partial-sum ascent-descent 
code as the LRMax code of a 321-avoiding permutation. 
Using the same Dyck path to illustrate, again 
$\A=(1,2,4,6,9),\ \D=(1,2,3,6,9)$. With $(\A,\D)$ as LRMax 
code, the left-to-right-maxima are given by \A with $n$ appended, 
their positions by $\D+1$ with 1 prepended.
Thus we have the partial permutation
\[
\begin{array}{cccccccccc}
 \textrm{{\scriptsize 1}} & \textrm{{\scriptsize 2}} & \textrm{{\scriptsize 3}} & 
    \textrm{{\scriptsize 4}}& \textrm{{\scriptsize 5}} & \textrm{{\scriptsize 
    6}} & 
    \textrm{{\scriptsize 7}} & \textrm{{\scriptsize 8}} 
    & \textrm{{\scriptsize 9}} & \textrm{{\scriptsize 10}}   \\[-0.8ex]
   1 & 2 & 4 & 6 &   &  & 9 &   &  & 10  
\end{array}
\]
and filling in the missing entries in increasing order gives the 
image permutation: \linebreak
$1\ 2\ 4\ 6\ 3\ 5\ 9\ 7\ 8\ 10$.

The Mansour-Deng-Du bijection $M$ is a bit more complicated and here we 
attempt to simplify its description, referring the reader to 
\cite{mansour2006} for full details of the original description. First 
label the upsteps of the Dyck path left to right and record the label 
on the first $u$ of each $uu$. 
Do likewise for the downsteps. For our running example (\ref{dyck1}), as already
noted in \S 5, the result is $(4,5,7,8)$ for the 
downsteps and $(3,5,7,8)$ for the upsteps and this pair forms the 
partial-sum ascent-descent code for $L(P)$. In \cite{mansour2006} 
this pair is denoted $(h,t)$ and is obtained by a different but equivalent process: a graphical 
construction based on the so-called $(x+y)$-\emph{labelling} of a Dyck path. 
Next, \cite{mansour2006} defines
$\si_{i}:=s_{h_{i}}s_{h_{i}-1}s_{h_{i}-2}\ldots 
s_{t_{i}}$  where $s_{j}$ is the transposition that interchanges $j$ and 
$j+1$, and goes on to form the image permutation as
\[
(1,2,\ \ldots,n) \si_{1}\si_{2}\ldots \si_{r},
\]
where $r$ is the length of $h$ (and $t$) and operations are performed 
left to right. The effect of these operations is simply to displace 
$h_{i}+1$ to the left in the list $(1,2,\ \ldots,n)$ so that it is in position 
$t_{i}$, this for each $i$ 
while leaving all other entries in the same relative order. A little 
thought shows that this is equivalent to using $h$ and $t$ as the excedance 
code to produce the image permutation. The example thus yields 
excedance values $h+1=(5,6,8,9)$ and excedance locations 
$t=(3,5,7,8)$, and 
so the image permutation is $1\ 2\ 5\ 3\ 6\ 4\ 8\ 9\ 7\ 10$.

\vspace*{8mm}

{\Large \textbf{8 \quad The identities $\mbf{M=B} \circ \mbf{L=K} \circ 
\mbf{L'}$  }  } It is now clear that $M=B \circ L$ because, as we have 
just 
seen, for a Dyck path $P$, the excedance code of $M(P)$ is the 
partial-sum ascent-descent code of $L(P)$ and the bijection $B$ uses the latter code
as an excedance code. To see that $M=K \circ L'$, equivalently, 
$K^{-1}\circ B = L' \circ L$, requires a little 
more work.

From the descriptions of $K$ and $B$ in the preceding section, we see 
that the following 4-step process transforms the LRMax code of a 321-avoiding 
permutation to its excedance code (writing the codes as 2-row 
matrices with the larger row on top):
\begin{enumerate}
    \item  append $n$  to the top row and and prepend 0 to the bottom 
    row 
    
    \item  add 1 to each entry of the bottom row

    \item  delete columns with same top and bottom entry

    \item  subtract 1 from each entry of the top row.
\end{enumerate}

For example, with $n=13$, 
\[
\left(
\begin{matrix}
    2 & 3 & 4 & 8 & 9 & 12 \\
    1 & 3 & 4 & 6 & 7 & 10
\end{matrix}
\right) 
\overset{(1)}{\longrightarrow}
\left(
\begin{matrix}
    2 & 3 & 4 & 8 & 9 & 12 & 13\\
    0 & 1 & 3 & 4 & 6 & 7 & 10
\end{matrix}
\right)
\overset{(2)}{\longrightarrow}
\left(
\begin{matrix}
    2 & 3 & 4 & 8 & 9 & 12 & 13\\
    1 & 2 & 4 & 5 & 7 & 8 & 11
\end{matrix}
\right)
\] 
\vspace*{-4mm} 
\begin{equation}
\overset{(3)}{\longrightarrow}
\left(
\begin{matrix}
    2 & 3  & 8 & 9 & 12 & 13\\
    1 & 2  & 5 & 7 & 8 & 11
\end{matrix}
\right)
\overset{(4)}{\longrightarrow}
\left(
\begin{matrix}
    1 & 2  & 7 & 8 & 11 & 12\\
    1 & 2  & 5 & 7 & 8  & 11
\end{matrix}
\right)
\label{kb}
\end{equation}
If $P$ is a 321-avoiding permutation and $D_{1},D_{2}$ are the Dyck 
paths corresponding to its LRMax and excedance codes respectively, 
then $K(D_{1})=B(D_{2})$ and so $D_{2}=B^{-1}\circ K(D_{1})$. We wish 
to trace the effects of the these 4 steps on the Dyck path $D_{1}$ and 
show that they produce $L \circ L' (D_{1})$; we can then conclude that 
$B^{-1}\circ K=L \circ L' $ or, taking inverses, that $K^{-1}\circ B=L' 
\circ L$, as desired.

The trick is to translate Dyck paths to path pairs using $\phi$. The 
composite bijection ``partial-sum ascent-descent code $\rightarrow$ Dyck 
path $\rightarrow$ path pair'' has a simple description as illustrated 
for the first entry in (\ref{kb}), with $n=13$:
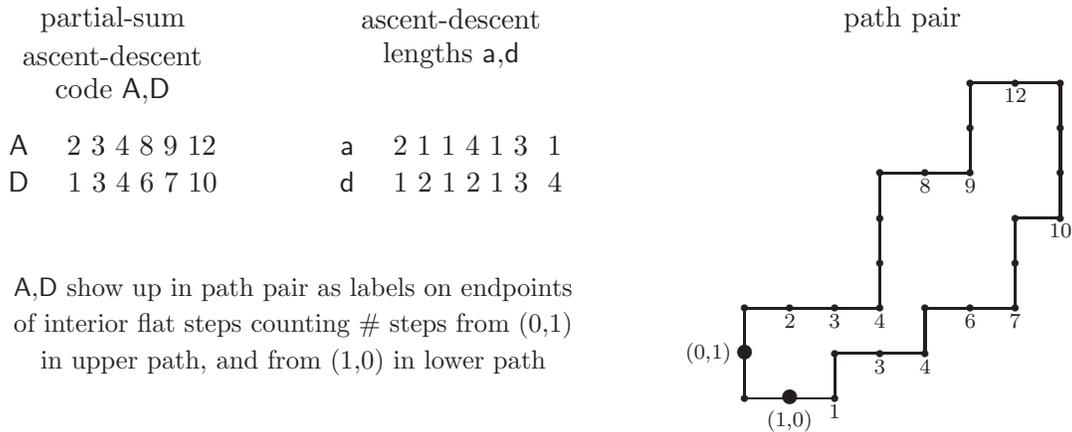
\vspace{1mm}

\Einheit=0.6cm
\[
\Pfad(5,0),22111222112211\endPfad
\Pfad(5,0),11211211221222\endPfad
\DuennPunkt(5,0)
\NormalPunkt(5,1)
\DuennPunkt(5,2)
\NormalPunkt(6,0)
\DuennPunkt(6,2)
\DuennPunkt(7,0)
\DuennPunkt(7,1)
\DuennPunkt(7,2)
\DuennPunkt(8,1)
\DuennPunkt(8,2)
\DuennPunkt(8,3)
\DuennPunkt(8,4)
\DuennPunkt(8,5)
\DuennPunkt(9,2)
\DuennPunkt(9,1)
\DuennPunkt(9,5)
\DuennPunkt(10,2)
\DuennPunkt(10,5)
\DuennPunkt(10,6)
\DuennPunkt(10,7)
\DuennPunkt(11,3)
\DuennPunkt(11,2)
\DuennPunkt(11,4)
\DuennPunkt(11,7)
\DuennPunkt(12,5)
\DuennPunkt(12,4)
\DuennPunkt(12,6)
\DuennPunkt(12,7)
\Label\l{{\textrm{\scriptsize (0,1)}}}(4.7,1)
\Label\u{{\textrm{\scriptsize (1,0)}}}(6,0)
\Label\o{{\textrm{\scriptsize 2}}}(6,1.3)
\Label\o{{\textrm{\scriptsize 3}}}(7,1.3)
\Label\o{{\textrm{\scriptsize 4}}}(8,1.3)
\Label\o{{\textrm{\scriptsize 8}}}(9,4.3)
\Label\o{{\textrm{\scriptsize 9}}}(10,4.3)
\Label\o{{\textrm{\scriptsize 12}}}(11,6.3)
\Label\o{{\textrm{\scriptsize 1}}}(7,-0.7)
\Label\o{{\textrm{\scriptsize 3}}}(8,0.3)
\Label\o{{\textrm{\scriptsize 4}}}(9,0.3)
\Label\o{{\textrm{\scriptsize 6}}}(10,1.3)
\Label\o{{\textrm{\scriptsize 7}}}(11,1.3)
\Label\o{{\textrm{\scriptsize 10}}}(12,3.3)
\Label\o{{\textrm{\small \A\quad 2\ 3\  4\  8\  9\  12}}}(-9,5.2)
\Label\o{{\textrm{\small \D\quad 1\ 3\  4\  6\  7\  10}}}(-9,4.4)
\Label\o{{\textrm{\small partial-sum}}}(-9,8)
\Label\o{{\textrm{\small ascent-descent}}}(-9,7.2)
\Label\o{{\textrm{\small code \A,\D}}}(-9,6.4)
\Label\o{{\textrm{\small ascent-descent}}}(-1.5,8)
\Label\o{{\textrm{\small lengths \a,\d}}}(-1.5,7.2)
\Label\o{{\textrm{\small path pair}}}(8.5,8)
\Label\o{{\textrm{\small \a\quad 2\ 1\  1\  4\  1\  3 \ 1}}}(-1.5,5.2)
\Label\o{{\textrm{\small \d\quad 1\ 2\  1\  2\  1\  3 \ 4}}}(-1.5,4.4)
\Label\o{{\textrm{\footnotesize \A,\D show up in path pair as labels on 
endpoints}}}(-5,2)
\Label\o{{\textrm{\footnotesize of interior flat steps counting \# steps from 
(0,1)}}}(-5,1.2)
\Label\o{{\textrm{\footnotesize  in upper path, and from (1,0) in lower 
path}}}(-5,0.4)
\]

Now we can describe the effect of the 4-step process on the path pair:
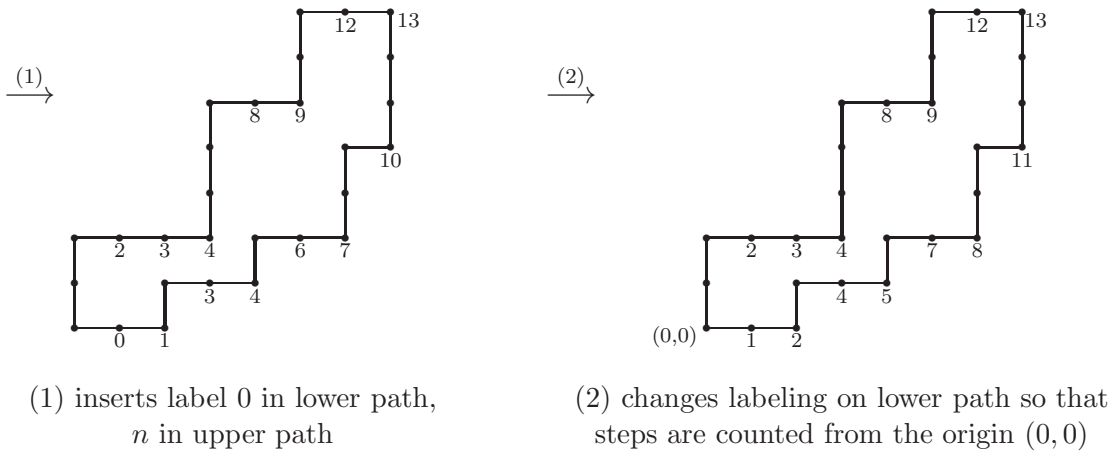
\Einheit=0.6cm
\[
\Label\o{\overset{(1)}{\longrightarrow}}(-11,5)
\Label\o{\overset{(2)}{\longrightarrow}}(1,5)
\Pfad(-10,0),22111222112211\endPfad
\Pfad(-10,0),11211211221222\endPfad
\DuennPunkt(-10,0)
\DuennPunkt(-10,1)
\DuennPunkt(-10,2)
\DuennPunkt(-9,0)
\DuennPunkt(-9,2)
\DuennPunkt(-8,0)
\DuennPunkt(-8,1)
\DuennPunkt(-8,2)
\DuennPunkt(-7,1)
\DuennPunkt(-7,2)
\DuennPunkt(-7,3)
\DuennPunkt(-7,4)
\DuennPunkt(-7,5)
\DuennPunkt(-6,2)
\DuennPunkt(-6,1)
\DuennPunkt(-6,5)
\DuennPunkt(-5,2)
\DuennPunkt(-5,5)
\DuennPunkt(-5,6)
\DuennPunkt(-5,7)
\DuennPunkt(-4,3)
\DuennPunkt(-4,2)
\DuennPunkt(-4,4)
\DuennPunkt(-4,7)
\DuennPunkt(-3,5)
\DuennPunkt(-3,4)
\DuennPunkt(-3,6)
\DuennPunkt(-3,7)
\Label\o{{\textrm{\scriptsize 2}}}(-9,1.3)
\Label\o{{\textrm{\scriptsize 3}}}(-8,1.3)
\Label\o{{\textrm{\scriptsize 4}}}(-7,1.3)
\Label\o{{\textrm{\scriptsize 8}}}(-6,4.3)
\Label\o{{\textrm{\scriptsize 9}}}(-5,4.3)
\Label\o{{\textrm{\scriptsize 12}}}(-4,6.3)
\Label\o{{\textrm{\scriptsize 1}}}(-8,-0.7)
\Label\o{{\textrm{\scriptsize 3}}}(-7,0.3)
\Label\o{{\textrm{\scriptsize 4}}}(-6,0.3)
\Label\o{{\textrm{\scriptsize 6}}}(-5,1.3)
\Label\o{{\textrm{\scriptsize 7}}}(-4,1.3)
\Label\o{{\textrm{\scriptsize 10}}}(-3,3.3)
\Label\o{{\textrm{\scriptsize 0}}}(-9,-0.7)
\Label\o{{\textrm{\scriptsize 13}}}(-2.6,6.4)
\Label\o{{\textrm{\small  (1) inserts label 0 in lower path,}}}(-6.5,-2)
\Label\o{{\textrm{\small $n$ in upper path}}}(-6.5,-2.8)
\Pfad(4,0),22111222112211\endPfad
\Pfad(4,0),11211211221222\endPfad
\DuennPunkt(4,0)
\DuennPunkt(4,1)
\DuennPunkt(4,2)
\DuennPunkt(5,0)
\DuennPunkt(5,2)
\DuennPunkt(6,0)
\DuennPunkt(6,1)
\DuennPunkt(6,2)
\DuennPunkt(7,1)
\DuennPunkt(7,2)
\DuennPunkt(7,3)
\DuennPunkt(7,4)
\DuennPunkt(7,5)
\DuennPunkt(8,2)
\DuennPunkt(8,1)
\DuennPunkt(8,5)
\DuennPunkt(9,2)
\DuennPunkt(9,5)
\DuennPunkt(9,6)
\DuennPunkt(9,7)
\DuennPunkt(10,3)
\DuennPunkt(10,2)
\DuennPunkt(10,4)
\DuennPunkt(10,7)
\DuennPunkt(11,5)
\DuennPunkt(11,4)
\DuennPunkt(11,6)
\DuennPunkt(11,7)
\Label\l{{\textrm{\scriptsize (0,0)}}}(3.8,-0.2)
\Label\o{{\textrm{\scriptsize 2}}}(5,1.3)
\Label\o{{\textrm{\scriptsize 3}}}(6,1.3)
\Label\o{{\textrm{\scriptsize 4}}}(7,1.3)
\Label\o{{\textrm{\scriptsize 8}}}(8,4.3)
\Label\o{{\textrm{\scriptsize 9}}}(9,4.3)
\Label\o{{\textrm{\scriptsize 12}}}(10,6.3)
\Label\o{{\textrm{\scriptsize 13}}}(11.3,6.4)
\Label\o{{\textrm{\scriptsize 1}}}(5,-0.7)
\Label\o{{\textrm{\scriptsize 2}}}(6,-0.7)
\Label\o{{\textrm{\scriptsize 4}}}(7,0.3)
\Label\o{{\textrm{\scriptsize 5}}}(8,0.3)
\Label\o{{\textrm{\scriptsize 7}}}(9,1.3)
\Label\o{{\textrm{\scriptsize 8}}}(10,1.3)
\Label\o{{\textrm{\scriptsize 11}}}(11,3.3)
\Label\o{{\textrm{\small  (2) changes labeling on lower path so that}}}(7,-2)
\Label\o{{\textrm{\small steps are counted from the origin $(0,0)$}}}(7,-2.8)
\]

\vspace*{1mm}
\[
\Label\o{\overset{(3)}{\longrightarrow}}(-11,5)
\Label\o{\overset{(4)}{\longrightarrow}}(1,5)
\Pfad(-10,0),22112222112211\endPfad
\Pfad(-10,0),11221211221222\endPfad
\DuennPunkt(-10,0)
\DuennPunkt(-10,1)
\DuennPunkt(-10,2)
\DuennPunkt(-9,0)
\DuennPunkt(-9,2)
\DuennPunkt(-8,0)
\DuennPunkt(-8,1)
\DuennPunkt(-8,2)
\DuennPunkt(-8,3)
\DuennPunkt(-8,4)
\DuennPunkt(-8,5)
\DuennPunkt(-8,6)
\DuennPunkt(-7,2)
\DuennPunkt(-7,3)
\DuennPunkt(-7,6)
\DuennPunkt(-6,3)
\DuennPunkt(-6,6)
\DuennPunkt(-6,7)
\DuennPunkt(-6,8)
\DuennPunkt(-5,3)
\DuennPunkt(-5,4)
\DuennPunkt(-5,5)
\DuennPunkt(-5,8)
\DuennPunkt(-4,5)
\DuennPunkt(-4,6)
\DuennPunkt(-4,7)
\DuennPunkt(-4,8)
\Label\o{{\textrm{\scriptsize 2}}}(-9,1.3)
\Label\o{{\textrm{\scriptsize 3}}}(-8.3,1.3)
\Label\o{{\textrm{\scriptsize 8}}}(-7,5.3)
\Label\o{{\textrm{\scriptsize 9}}}(-6,5.3)
\Label\o{{\textrm{\scriptsize 12}}}(-5,7.3)
\Label\o{{\textrm{\scriptsize 1}}}(-9,-0.7)
\Label\o{{\textrm{\scriptsize 2}}}(-8,-0.7)
\Label\o{{\textrm{\scriptsize 5}}}(-7,1.3)
\Label\o{{\textrm{\scriptsize 7}}}(-6,2.3)
\Label\o{{\textrm{\scriptsize 8}}}(-5,2.3)
\Label\o{{\textrm{\scriptsize 11}}}(-4,4.3)
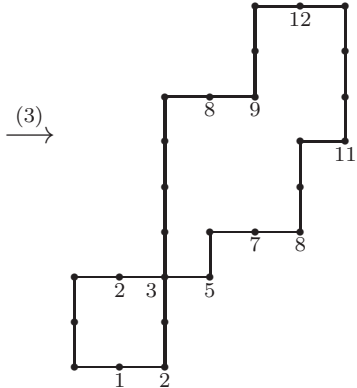
\Label\o{{\textrm{\small (3) swings flat steps forming}}}(-7.5,-2)
\Label\o{{\textrm{\small sides of unit square }}}(-7.5,-2.8)
\Label\o{{\textrm{\small to vertical steps}}}(-7.5,-3.6)
\Pfad(4,0),21122221122111\endPfad
\Pfad(4,0),11122121122122\endPfad
\DuennPunkt(4,0)
\DuennPunkt(4,1)
\DuennPunkt(5,0)
\DuennPunkt(5,1)
\DuennPunkt(6,0)
\DuennPunkt(6,1)
\DuennPunkt(6,2)
\DuennPunkt(6,3)
\DuennPunkt(6,4)
\DuennPunkt(6,5)
\DuennPunkt(7,0)
\DuennPunkt(7,1)
\DuennPunkt(7,2)
\DuennPunkt(7,5)
\DuennPunkt(8,2)
\DuennPunkt(8,3)
\DuennPunkt(8,5)
\DuennPunkt(8,6)
\DuennPunkt(8,7)
\DuennPunkt(9,3)
\DuennPunkt(9,7)
\DuennPunkt(10,3)
\DuennPunkt(10,4)
\DuennPunkt(10,5)
\DuennPunkt(10,7)
\DuennPunkt(11,5)
\DuennPunkt(11,6)
\DuennPunkt(11,7)
\Label\o{{\textrm{\scriptsize 1}}}(5,0.3)
\Label\o{{\textrm{\scriptsize 2}}}(6,0.3)
\Label\o{{\textrm{\scriptsize 7}}}(7,4.3)
\Label\o{{\textrm{\scriptsize 8}}}(8,4.3)
\Label\o{{\textrm{\scriptsize 11}}}(9,6.3)
\Label\o{{\textrm{\scriptsize 12}}}(10,6.3)
\Label\o{{\textrm{\scriptsize 1}}}(6,-0.7)
\Label\o{{\textrm{\scriptsize 2}}}(7,-0.7)
\Label\o{{\textrm{\scriptsize 5}}}(8,1.3)
\Label\o{{\textrm{\scriptsize 7}}}(9,2.3)
\Label\o{{\textrm{\scriptsize 8}}}(10,2.3)
\Label\o{{\textrm{\scriptsize 11}}}(11,4.3)
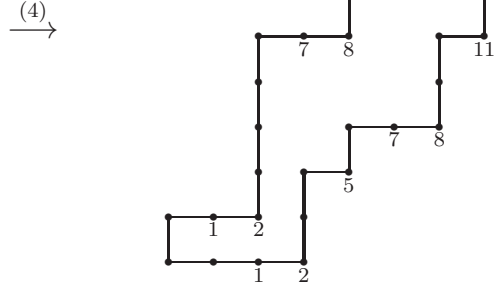
\Label\o{{\textrm{\small (4) subtracts 1 from each label in upper path,}}}(6,-2)
\Label\o{{\textrm{\small and to restore the counting of steps from}}}(6.5,-2.8)
\Label\o{{\textrm{\small (0,1) in the upper path and from (1,0) in}}}(6.5,-3.6)
\Label\o{{\textrm{\small the lower path amounts to rotating the initial}}}(6.5,-4.4)
\Label\o{{\textrm{\small (vertical) step in each path $90^{\circ}$ 
counterclockwise}}}(6.5,-5.2)
\]
\vspace*{2mm}

\noindent It is evident that the final result is indeed $L\circ L'$ applied to 
the initial path pair because the initial path pair is obtained from 
the final path pair by applying $L'\circ L$ as described in \S 6. 
Thus we have shown that $K^{-1} \circ B =  L' \circ L $.
\vspace*{8mm}

{\Large \textbf{9 \quad Analysis of $\mbf{K^{-1}} \circ \mbf{B} =  \mbf{L'} \circ \mbf{L} $.  }  }

Recall that $R$ is the ``reverse path'' involution on Dyck paths and 
is also, under $\phi$, the 
``rotate $180^{\circ}$'' involution on path pairs.  
\begin{theorem}
 On Dyck $n$-paths,
 \[
 (L' \circ L)^{n}=R.
 \]
\end{theorem}
\begin{cor}
   For $n\ge 3$, the permutation $L' \circ L$ on  Dyck $n$-paths has order $2n$.
\end{cor}
\textbf{Proof of Corollary}\quad Since $R$ is an involution, the theorem shows that the order of $L' 
\circ L$ divides $2n$. The assertion can be checked directly for $n=3$ 
and for $n\ge 4$, the orbit of the Dyck path $u^{n-1}d^{n-1}ud$ 
(exponents denote repetition) has size $2n$. \qed

\textbf{Proof of Theorem}\quad 
We will consider the effect of repeated application of $L' \circ L $ 
on a path pair $P$ of size $n+1$.
Recall from \S 6 that $L'\circ L(P)$ is obtained as follows:

\noindent (i) identify the last (northeasternmost) step in the upper 
path and the southwesternmost step in the lower path.

\noindent (ii) identify each pair of vertical steps that form opposite 
sides of a unit square.

Then change the two steps in (i) from flat to vertical and all steps 
in (ii) (if any) from vertical to flat.

Now consider a path pair as a linkage composed of rods of 
unit length that must always be aligned either flat or vertical, 
hinged at the vertices. Applying $L'\circ L$ then simply changes the 
alignment of some of the rods (steps) but preserves their identity, 
that is, one may track the progress of a particular step or vertex 
under repeated applications of $L'\circ L$.

Let us count steps in a path pair clockwise from the origin. Thus the 
first and $(n+2)$nd steps initiate the upper and lower paths 
respectively and both are necessarily vertical. If a step $S$ is the 
$i$th step in a path pair $P$, then $S$ becomes step number $i+1$ 
(mod $2n+2$) in $L'\circ L(P)$. In particular, under $(L'\circ 
L)^{n}$, the initial steps in the upper and lower paths become their 
terminal steps respectively while every other step passes from its original path 
to the other one. When it does so, we will say it ``turns the corner''.

It is clear that, under repeated applications of $L'\circ L$, each vertical 
step must get flattened before it turns the corner and, once flat, a 
step stays flat until it turns the corner (when, of course, it 
becomes vertical). So the crux of the matter is whether or not a step 
gets flattened \emph{after} it turns the corner. To show that the 
effect of $(L'\circ L)^{n}$ is to rotate a path pair $180^{\circ}$, 
we must show 
\begin{prop}
   Let $P$ be a path pair of size $n+1$.  Under $n$ applications of  $L'\circ L$, a step in $P$ gets
    flattened after it turns the corner if and only if it  is immediately 
    preceded by a flat step in the original path pair.
\end{prop}
\textbf{Proof}\quad
A \emph{minimal diagonal} in a path pair is a line segment of slope 
1 ($45^{\circ}$) joining two distinct vertices (either in the same or 
different paths) and lying strictly inside the path pair except at 
its endpoints. A simple count shows that there are exactly $n$ minimal 
diagonals in a path pair of size $n+1$. Given a 
minimal diagonal, let $V_{1}<V_{2}$ (counting clockwise) denote its 
endpoints and let $S_{1},S_{2}$ denote the steps initiated 
(clockwise) by $V_{1},V_{2}$ respectively. The key observation is 
that as $V_{1},V_{2}$ progress under repeated applications of $L'\circ 
L$, they remain endpoints of a minimal diagonal until 
$S_{1},S_{2}$ form the vertical sides of a unit square. As the next 
paragraph shows, this will 
always happen at $(L'\circ L)^{i}(P)$ for some $i,\ 0\le i \le n-1$. 
The next application of $L'\circ L$ will then flatten both $S_{1}$ and 
$S_{2}$. Furthermore, by tracing backwards, every instance of a pair of steps forming 
vertical sides of a unit square in the set $\{(L'\circ 
L)^{i}(P)\}_{i=0}^{n-1}$  arises in this way from a minimal diagonal 
of $P$.

Applying $L'\circ L$ changes the length of a minimal diagonal by at 
most 1. Specifically, if $V_{1},V_{2}$ are interior points of 
different paths and $\overrightarrow{V_{1}V_{2}}$ points southwest, 
the length increases by 1; if $V_{1},V_{2}$ are in the same path and 
$V_{2}$ is not the path's terminal point, the length stays the same; 
otherwise, the length decreases by 1. It follows that a minimal 
diagonal can survive at most $n-1$ applications of $L'\circ L$ before being 
destroyed at the next application. 

If (initially) $S_{1}$ lies in the upper path and $S_{2}$ in the lower 
path, then both are vertical and get flattened before 
either turns the corner. If $S_{1}$ lies in the lower path and $S_{2}$ in the 
upper, then each is preceded by a flat step and flattening occurs 
after both $S_{1}$ and $S_{2}$ have turned the corner. If 
$S_{1},\,S_{2}$ both lie in the same path (either upper or lower), 
then $S_{1}$ is vertical, $S_{2}$ is preceded by a flat step and flattening 
occurs before $S_{1}$ turns the corner and after $S_{2}$ does so. The 
Proposition follows.

An example with $n=8$ is shown along with the progress of 3 of the 8 
minimal diagonals, using a different color for each one.

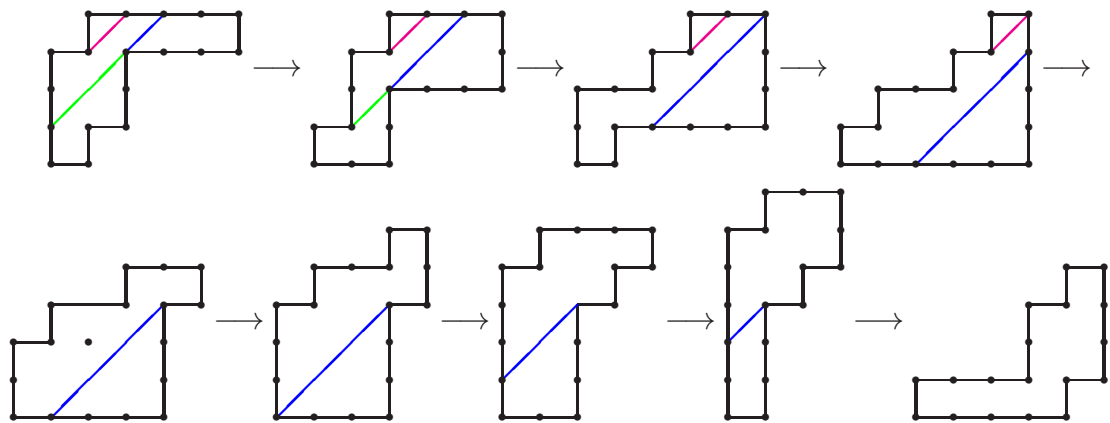
\vspace*{-2mm}

\Einheit=0.5cm
\[
\Pfad(-13,0),222121111\endPfad
\Pfad(-13,0),121221112\endPfad
\Pfad(-6,0),212212111\endPfad
\Pfad(-6,0),112211122\endPfad
\Pfad(1,0),221121211\endPfad
\Pfad(1,0),121111222\endPfad
\Pfad(8,0),212112121\endPfad
\Pfad(8,0),111112222\endPfad
\magenta{
\Pfad(-12,3),3\endPfad
\Pfad(-4,3),3\endPfad
\Pfad(4,3),3\endPfad
\Pfad(12,3),3\endPfad }
\green{
\Pfad(-13,1),33\endPfad
\Pfad(-5,1),3\endPfad }
\blue{
\Pfad(-11,3),3\endPfad
\Pfad(-4,2),33\endPfad
\Pfad(3,1),333\endPfad
\Pfad(10,0),333\endPfad }
\DuennPunkt(-13,0)
\DuennPunkt(-13,1)
\DuennPunkt(-13,2)
\DuennPunkt(-13,3)
\DuennPunkt(-12,0)
\DuennPunkt(-12,1)
\DuennPunkt(-12,3)
\DuennPunkt(-12,4)
\DuennPunkt(-11,1)
\DuennPunkt(-11,2)
\DuennPunkt(-11,3)
\DuennPunkt(-11,4)
\DuennPunkt(-10,3)
\DuennPunkt(-10,4)
\DuennPunkt(-9,3)
\DuennPunkt(-9,4)
\DuennPunkt(-8,3)
\DuennPunkt(-8,4)
\DuennPunkt(-6,1)
\DuennPunkt(-6,0)
\DuennPunkt(-5,0)
\DuennPunkt(-5,1)
\DuennPunkt(-5,2)
\DuennPunkt(-5,3)
\DuennPunkt(-4,0)
\DuennPunkt(-4,1)
\DuennPunkt(-4,2)
\DuennPunkt(-4,3)
\DuennPunkt(-4,4)
\DuennPunkt(-3,2)
\DuennPunkt(-3,4)
\DuennPunkt(-2,2)
\DuennPunkt(-2,4)
\DuennPunkt(-1,2)
\DuennPunkt(-1,3)
\DuennPunkt(-1,4)
\Label\o{\longrightarrow}(-7,2)
\Label\o{\longrightarrow}(0,2)
\Label\o{\longrightarrow}(7,2)
\Label\o{\longrightarrow}(14,2)
\DuennPunkt(1,0)
\DuennPunkt(1,1)
\DuennPunkt(1,2)
\DuennPunkt(2,0)
\DuennPunkt(2,1)
\DuennPunkt(2,2)
\DuennPunkt(3,1)
\DuennPunkt(3,2)
\DuennPunkt(3,3)
\DuennPunkt(4,1)
\DuennPunkt(4,3)
\DuennPunkt(4,4)
\DuennPunkt(5,1)
\DuennPunkt(5,4)
\DuennPunkt(6,1)
\DuennPunkt(6,2)
\DuennPunkt(6,3)
\DuennPunkt(6,4)
\DuennPunkt(8,0)
\DuennPunkt(8,1)
\DuennPunkt(9,0)
\DuennPunkt(9,1)
\DuennPunkt(9,2)
\DuennPunkt(10,0)
\DuennPunkt(10,2)
\DuennPunkt(11,0)
\DuennPunkt(11,2)
\DuennPunkt(11,3)
\DuennPunkt(12,0)
\DuennPunkt(12,3)
\DuennPunkt(12,4)
\DuennPunkt(13,0)
\DuennPunkt(13,1)
\DuennPunkt(13,2)
\DuennPunkt(13,3)
\DuennPunkt(13,4)
\]
\Einheit=0.5cm
\[
\Pfad(-14,0),221211211\endPfad
\Pfad(-14,0),111122212\endPfad
\Pfad(-7,0),222121121\endPfad
\Pfad(-7,0),111222122\endPfad
\Pfad(-1,0),222212111\endPfad
\Pfad(-1,0),112221212\endPfad
\Pfad(5,0),222221211\endPfad
\Pfad(5,0),122212122\endPfad
\Pfad(10,0),211122121\endPfad
\Pfad(10,0),111121222\endPfad
\blue{
\Pfad(-13,0),333\endPfad
\Pfad(-7,0),333\endPfad
\Pfad(-1,1),33\endPfad
\Pfad(5,2),3\endPfad }
\DuennPunkt(-14,0)
\DuennPunkt(-14,1)
\DuennPunkt(-14,2)
\DuennPunkt(-13,0)
\DuennPunkt(-13,2)
\DuennPunkt(-13,3)
\DuennPunkt(-12,0)
\DuennPunkt(-12,2)
\DuennPunkt(-11,0)
\DuennPunkt(-11,3)
\DuennPunkt(-11,4)
\DuennPunkt(-10,0)
\DuennPunkt(-10,1)
\DuennPunkt(-10,3)
\DuennPunkt(-10,4)
\DuennPunkt(-10,2)
\DuennPunkt(-9,3)
\DuennPunkt(-9,4)
\DuennPunkt(-7,0)
\DuennPunkt(-7,1)
\DuennPunkt(-7,2)
\DuennPunkt(-7,3)
\DuennPunkt(-6,0)
\DuennPunkt(-6,3)
\DuennPunkt(-6,4)
\DuennPunkt(-5,0)
\DuennPunkt(-5,4)
\DuennPunkt(-4,0)
\DuennPunkt(-4,1)
\DuennPunkt(-4,2)
\DuennPunkt(-4,3)
\DuennPunkt(-4,4)
\DuennPunkt(-4,5)
\DuennPunkt(-3,3)
\DuennPunkt(-3,4)
\DuennPunkt(-3,5)
\DuennPunkt(-1,0)
\DuennPunkt(-1,1)
\DuennPunkt(-1,2)
\DuennPunkt(-1,3)
\DuennPunkt(-1,4)
\DuennPunkt(0,0)
\DuennPunkt(0,4)
\DuennPunkt(0,5)
\Label\o{\longrightarrow}(-8,2)
\Label\o{\longrightarrow}(-2,2)
\Label\o{\longrightarrow}(4,2)
\Label\o{\longrightarrow}(9,2)
\DuennPunkt(1,0)
\DuennPunkt(1,1)
\DuennPunkt(1,2)
\DuennPunkt(1,5)
\DuennPunkt(2,3)
\DuennPunkt(2,4)
\DuennPunkt(2,5)
\DuennPunkt(3,4)
\DuennPunkt(3,5)
\DuennPunkt(5,0)
\DuennPunkt(5,1)
\DuennPunkt(5,2)
\DuennPunkt(5,3)
\DuennPunkt(5,4)
\DuennPunkt(5,5)
\DuennPunkt(6,0)
\DuennPunkt(6,1)
\DuennPunkt(6,2)
\DuennPunkt(6,3)
\DuennPunkt(6,5)
\DuennPunkt(6,6)
\DuennPunkt(7,3)
\DuennPunkt(7,4)
\DuennPunkt(7,6)
\DuennPunkt(8,4)
\DuennPunkt(8,5)
\DuennPunkt(8,6)
\DuennPunkt(10,0)
\DuennPunkt(10,1)
\DuennPunkt(11,0)
\DuennPunkt(11,1)
\DuennPunkt(12,0)
\DuennPunkt(12,1)
\DuennPunkt(13,0)
\DuennPunkt(13,1)
\DuennPunkt(13,2)
\DuennPunkt(13,3)
\DuennPunkt(14,0)
\DuennPunkt(14,1)
\DuennPunkt(14,3)
\DuennPunkt(14,4)
\DuennPunkt(15,1)
\DuennPunkt(15,2)
\DuennPunkt(15,3)
\DuennPunkt(15,4)
\]

That $K^{-1} \circ B$ turns out to be a product of two ``nice'' 
involutions may be somewhat unexpected but recall that every 
(ordinary)
permutation can be expressed as a product of two involutions \cite{2invols}.

\end{document}